\documentclass{amsart}
\usepackage{amsmath,  amssymb,  graphicx, tikz, booktabs}
\usetikzlibrary{arrows.meta, cd}

\usepackage{hyperref}

\hypersetup{hypertexnames=false}


\newcommand{\PP}[1]{\textsc{Subsets}#1}
\newcommand{\PPs}[1]{\textsc{PSubsets}#1}

\usepackage{pifont}

\newcommand{\ZZ}{\mathbb{Z}}
\newcommand{\RR}{\mathbb{R}}

\newcommand{\GL}[1]{\mathit{GL}_{#1}}

\newcommand{\SU}[1]{\mathit{SU(#1)}}

\newcommand{\lr}[3]{c_{#1;#2}^{#3}}   % LR coefficiients
\newcommand{\tm}[3]{c(#1; #2; #3)} % triple multiplicities
\newcommand{\mult}[2]{\mathit{mult}\left(#1; #2\right)} % multiplicity of #1 in #2

\newcommand{\lspace}[1]{\mathcal{L}_{\scriptscriptstyle{\mathrm{#1}}}}
\newcommand{\lattice}[1]{\Lambda_{\scriptscriptstyle{\mathrm{#1}}}}
\newcommand{\cone}[1]{{\textsf{#1}}}
\newcommand{\latof}[1]{\operatorname{lat}(#1)} % lattice points of a cone

\newcommand{\proj}{\mathit{pr}}
\newcommand{\pr}[1]{\proj_{\scriptscriptstyle{\mathrm{#1}}}}
\newcommand{\chc}[1]{\mathcal{K}_{\scriptscriptstyle{\mathrm{#1}}}}

\newcommand{\Horn}{\cone{H}}

\newcommand{\chKLR}{\chc{LRC}}

\newcommand{\TM}{\cone{TM}}
\newcommand{\latTM}{\latof{\TM}}
\newcommand{\chK}{\chc{TM}}

\newcommand{\BZ}{\cone{BZ}}
\newcommand{\coneBZ}{\cone{BZ}} % ST: I added this to fix notation issue from line 5.
\newcommand{\spanBZ}{\lspace{BZ}}
\newcommand{\latBZ}{\latof{\BZ}}

\newcommand{\Gg}{G_\textrm{g}}
\newcommand{\Gl}{G_\textrm{l}}
\newcommand{\Glg}{G_\textrm{lg}}
\newcommand{\upgroup}{G_{\textrm{\tiny{BZ}}}}

\newcommand{\ray}[1]{\protect\overrightarrow{#1}}
\newcommand{\upray}[1]{\Delta_{\ray{#1}}}
\newcommand{\lrray}[1]{\ray{#1^{\bullet}}}
\newcommand{\ulrray}{\ray{\lambda}}
\newcommand{\vlrray}{\ray{\mu}}
\newcommand{\upcone}{\widehat{C}}
\newcommand{\ltriangle}{\triangleleft}
\newcommand{\rtriangle}{\triangleright}

\renewcommand{\S}[1]{\mathfrak{S}_{#1}}

\newcommand{\Vol}[1]{\operatorname{Vol}_{#1}}% volume wrt a lattice
\newcommand{\sol}{\textsc{Sol}}% set of solutions  

  \tikzstyle{medv}=[draw,circle,fill=white,minimum size=10pt,inner sep=0pt]
  \tikzstyle{central vertex}=[draw,circle,fill=white,minimum size=20pt,inner sep=0pt]
\tikzstyle{w}=[draw, circle, fill=white,minimum size=6pt,inner sep=0pt]
\tikzstyle{b}=[draw, circle, fill=black,minimum size=6pt,inner sep=0pt]
\tikzstyle{zero}=[minimum size=6pt]
\tikzstyle{r}=[draw,circle,fill=white,minimum size=3pt,inner sep=0pt]
\tikzstyle{rb}=[draw,circle,fill=black,minimum size=1pt,inner sep=0pt]
\tikzstyle{z}=[draw,circle,fill=black,minimum size=1pt,inner sep=0pt]

\newcommand{\facediagram}{
  \tikzstyle{n}=[draw,circle,fill=black,minimum size=2pt,inner sep=0pt]
  \tikzstyle{g}=[draw,circle,fill=lightgray,minimum size=6pt,inner sep=0pt]
        \draw (0,0) node (Y3)   [Y3] {}
    -- ++(0:2.0cm) node (Z6) [Z6] {}
    -- ++(0:2.0cm) node (Z1) [Z1]  {}
    -- ++(0:2.0cm) node (Y1)  [Y1]  {};
    \draw (Y3)
    -- ++(60:2.0cm) node (Z5) [Z5]  {}
    -- ++(60:2.0cm) node (Z4) [Z4]   {}
    -- ++(60:2.0cm) node (Y2) [Y2]  {};
    \draw (Y2)
    -- ++(-60:2.0cm) node (Z3) [Z3]  {}
    -- ++(-60:2.0cm) node (Z2) [Z2] {};
    \draw (Z2) -- (Y1);
    \draw (Z5) -- (Z6);
    \draw (Z1) -- (Z2);
    \draw (Z3) -- (Z4);
    \draw (3, 1.73) node (X)  [X] {};
}

\newcommand{\BZT}[9]{
    \tikzstyle{r}=[fill=white,minimum size=10pt,inner sep=0pt]
     \path (0,0) node (Y3)  [r]  {$#9$}
    -- ++(0:2.0cm) node (Z6)    [r] {$#6$}
    -- ++(0:2.0cm) node (Z1)     [r] {$#1$}
    -- ++(0:2.0cm) node (Y1)   [r]   {$#7$};
    \path (Y3)
    -- ++(60:2.0cm) node (Z5)   [r]  {$#5$}
    -- ++(60:2.0cm) node (Z4)    [r]  {$#4$}
    -- ++(60:2.0cm) node (Y2)    [r] {$#8$};
    \path (Y2)
    -- ++(-60:2.0cm) node (Z3)   [r]  {$#3$}
    -- ++(-60:2.0cm) node (Z2)   [r] {$#2$};
    \draw (Z1) -- (Z2) -- (Z3) -- (Z4) -- (Z5) -- (Z6) -- (Z1) -- cycle;
    \draw (Z6) -- (Y3) -- (Z5);
    \draw (Z1) -- (Y1) -- (Z2);
    \draw (Z3) -- (Y2) -- (Z4);
 }

\newcommand{\fullBZT}[9]{
  \BZT{#1}{#2}{#3}{#4}{#5}{#6}{#7}{#8}{#9}
    \draw (Z6) -- (Y3) node [midway, scale=0.7, below] () {$\the\numexpr #6 + #9$};
    \draw (Z1) -- (Y1) node [midway, scale=0.7, below] () {$\the\numexpr #1 + #7$};
    \draw (Z2) -- (Y1) node [midway, scale=0.7, right] () {$\the\numexpr #2 + #7$};
    \draw (Z3) -- (Y2) node [midway, scale=0.7, right] () {$\the\numexpr #3 + #8$};
    \draw (Z4) -- (Y2) node [midway, scale=0.7, left] () {$\the\numexpr #4 + #8$};
    \draw (Z5) -- (Y3) node [midway, scale=0.7, left] () {$\the\numexpr #5 + #9$};
}

  \newcommand{\fdiagram}{
  \node (E) at (0,0)  [E] {};
  \node (C1) at (-30:1cm) [C1] {};
  \node (D1) at (30:1cm) [D1] {};
  \node (C2) at (90:1cm) [C2] {};
  \node (D3) at (150:1cm) [D3] {};
  \node (C3) at (210:1cm) [C3] {};
  \node (D5) at (270:1cm) [D5] {};
  \draw (E);
  \draw (C1)--(D1)--(C2)--(D3)--(C3)--(D5)-- (C1);
  }

  \newcommand{\lefttdiag}{
   \draw[very thick] (-60:0.86cm) -- (60:0.86cm) -- (180:0.86cm) -- cycle ;
  }

  \newcommand{\righttdiag}{
   \draw[very thick] (120:0.86cm) -- (240:0.86cm) -- (0:0.86cm) -- cycle ;
  }

\title{All linear symmetries of the $\SU{3}$ tensor multiplicities}
\author{Emmanuel Briand}\address{E. Briand and M.Rosas, Universidad de Sevilla}
\author{Mercedes Rosas}
\author{Stefan Trandafir}\address{S. Trandafir, Simon Fraser University, Vancouver}
\thanks{Partially supported by the Grant PID2020-117843GB-I00 funded by MCIN/AEI/10.13039/501100011033 and  by P20\_01056 (Junta de Andalucía/PAIDI 2020/FEDER)}
\date{\today}

\begin{document}

%\today
%\input{changes.txt} 
%\tableofcontents
%\cleardoublepage
%\setcounter{page}{1} 

\begin{abstract}
  The $\SU{3}$ tensor multiplicities are piecewise polynomial of degree $1$ in their labels. The pieces are the chambers of a complex of cones. We describe in detail this chamber complex and determine the group of all linear symmetries (of order $144$) for these tensor multiplicities. We represent the cells by diagrams showing clearly the inclusions  as well as the actions of the group of symmetries and of its remarkable subgroups.
\end{abstract}

\maketitle

%%%%%%%%%%%%%%%%%%%%%%%%%%%%%%%%%%%%%%%%%%%%%%%%%%%%%%%%%%%%%%%%%%%%%%%%%%%%%%%
%% Intro
%%----------------------------------------------------------------------

\section{Introduction}

The tensor multiplicities for $\SU{k}$ appear in several different contexts in nuclear and particle physics, as well as quantum interferometry, \cite{Elliott, Gellmann, Neeman, Mandeltsveig, ReckZeilingerBernstein, Wesslen}. For instance, in the classification of orbital states of particles in the nuclear shell model, the  $\SU{k}$ tensor multiplicities represent total angular momentum \cite{Elliott}. In quantum interferometry, the $\SU{k}$ tensor multiplicities appear in the decomposition of input states into a direct sum of input states in distinct, finite irreducible representations \cite{Sanders_et_al}. They also appear in elementary particle physics in \emph{the Eightfold way} pioneered by Gell-mann and Ne'eman \cite{Gellmann, Neeman}.

The $\SU{k}$ tensor multiplicities  also govern  the decomposition into irreducibles of  some restrictions of representations of  the symmetric group, describe the structure constants for the multiplication in the cohomology ring of the Grassmannian,  give a basis for the Schubert cycles, and are the structure constants for the multiplication of symmetric functions under the Schur basis, see \cite{Fulton}.

The $\SU{k}$ tensor multiplicities afford a number of symmetries. One of these symmetries is obvious: the multiplicity of $V_{n}$ in $V_{\ell} \otimes V_{m}$ is the same as in $V_{m} \otimes V_{\ell}$. Here $V_{\ell}$, $V_{m}$ \ldots are the irreducible representations of $\SU{k}$.
The $\SU{k}$ tensor multiplicities define a function on the lattice of the triples $(\ell, m, n)$ of weights for $\SU{k}$.
The aforementioned symmetry, $\ell \leftrightarrow m$, is \emph{linear} since it is given by a linear function on this lattice.
There is actually a well-known group of $12$ linear symmetries for the $\SU{k}$ tensor multiplicities (comprising $\ell \leftrightarrow m$), for any $k$.

In this paper, we focus on the case $k=3$.
Our main result is the determination, in an elementary way, of the full group of all linear symmetries for the $\SU{3}$ tensor multiplicities.
We find it has order $144$, surprisingly bigger than the group of $12$ \emph{general symmetries} that hold for $\SU{k}$ for any $k$.
We  determine the structure of this group:  it is   isomorphic to $\S{2} \times \left(\S{3} \wr \S{2}\right)$
%(equivalently $\S{2} \times \left((\S{3} \times \S{3}) \rtimes \S{2}\right)$). 
where $\S{n}$ are the symmetric groups.

It is known that, for any fixed $k$, the  $\SU{k}$ tensor multiplicities are given by piecewise polynomial formulas.
The pieces (domains of validity of the polynomial formulas) are the chambers (maximal cells) of a complex of polyhedral cones (the \emph{chamber complex}).
In the case $k=3$, we get a full, explicit description of the chamber complex from the combinatorial description (BZ triangles) given by Berenstein and Zelevinsky of these tensor multiplicities \cite{BerensteinZelevinsky}. We get  as well the action of the linear symmetries on the chamber complex.
We associate to each cell of the chamber complex a diagram, in such a way that the action of the linear symmetries on the cells can be read easily from their diagrams.

Following Berenstein and Zelevinsky, we actually work with the \emph{triple multiplicities of $\SU{3}$} rather than with the $\SU{3}$ tensor multiplicities themselves.  The triple multiplicities are an avatar of the tensor multiplicities, defined as the dimensions of the subspaces of invariants $(V_\ell \otimes V_m \otimes V_n)^{SU(3)}$. The advantage of this model is the straightforward description of the group of the $12$ general symmetries: is is generated by the 6 permutations of the factors $V_{\ell}$, $V_{m}$, $V_{n}$, and the duality involution.
Also the combinatorial description given by Berenstein and Zelevinsky is given in the setting of the triple multiplicities.

The tensor multiplicities for $\SU{k}$ are essentially the same as the tensor multiplicities for $GL_k$. We translate our results to the $GL_3$-setting, obtaining for the corresponding Littlewood-Richardson coefficients the full group of linear symmetries (it has order $288$) and the description of the chamber complex.

\subsubsection*{Relation with other works}

Two of us derived the full group of linear symmetries for the $\SU{3}$ tensor multiplicities in  the preprint \cite{BriandRosas}.
The proofs there involved some computer calculations.
This paper supersedes this preprint as it gets rid of all computer calculations.
It provides instead simple, clear proofs with only elementary calculations.

Later, \cite{CrampeEtAl} provided  a nice explanation for the existence of ``additional'' symmetries  in the $\SU{3}$ case (with respect to the general symmetries that hold for $\SU{k}$, for any $k$).

The first to provide complete explicit formulas for the  $\SU{3}$ tensor multiplicities is probably Mandel{'}tsve\u{\i}g \cite{Mandeltsveig}.
His description, however, lacks the clarity provided by the polyhedral geometric setting of the chamber complex, which puts order and structure to the formulas.  The chamber complex for the   $\SU{3}$ tensor multiplicities (actually for the Littlewood--Richardson coefficients of $\GL{3}$) was provided by Rassart \cite{Rassart}. Again, Rassart's derivation made use of some computer calculations. In our work, we get again the description of the chamber complex, but with no computer calculation involved, and relate its structure with the linear symmetries.

\subsubsection*{Structure of the paper}

Section \ref{sec:preliminaries} provides a presentation of the general framework.
Section \ref{sec: BZ triangles} introduces the BZ triangles for $\SU{3}$.
Section \ref{sec:symmetries} contains the main result of this work: the determination of all linear symmetries of the triple multiplicities of $\SU{3}$. 
In Section \ref{sec:chamber complex}, the full description of the chamber complex for the $\SU{3}$ tensor multiplicities is derived, as well as the action of the linear symmetries on it. 
Section \ref{sec:LR} is devoted to translating  the previous results  (symmetries and chamber complex) into the  $\GL{3}$ tensor multiplicities setting.
Section \ref{sec:volume} presents a determinantal formula for the $\SU{3}$ tensor multiplicities. Section \ref{sec:stability} uses the description of the $\SU{3}$ chamber complex to illustrate concretely known general stability properties of the Littlewood-Richardson coefficients.
Finally, Section \ref{sec:remarks} concludes with a  discussion of  some topics not contemplated  in the main bulk of the paper.

%%%%%%%%%%%%%%%%%%%%%%%%%%%%%%%%%%%%%%%%%%%%%%%%%%%%%%%%%%%%%%%%%%%%%%%%%%%%%%%
%% Preliminaries
%%----------------------------------------------------------------------

\section{Preliminaries}\label{sec:preliminaries}

\subsection{Representations of $\SU{3}$ and their labels}

The irreducible Lie group representations of $\SU{3}$ are naturally labelled by their highest weights. These are vectors in an abstract vector space (the dual of the complexified Lie algebra of $\SU{3}$). We will use as \emph{numeric} labels the vectors of coordinates $(\ell_1, \ell_2)$ of the highest weights in the basis of fundamental weights (\emph{Dynkin labels}).  We will denote with $V_{\ell}$ the irreducible representation of $\SU{3}$ with Dynkin label $\ell=(\ell_1, \ell_2)$. The Dynkin labels of the irreducible representations of $\SU{3}$ are the pairs of nonnegative integers.

The Dynkin label of the dual of the representation of $V_{\ell}$ with Dynkin label $\ell=(\ell_1, \ell_2)$ is $\ell$ read backwards,i.e.  $(\ell_2, \ell_1)$. We denote this by $\ell^*$.

\subsection{Triple multiplicities and tensor multiplicities}

Consider Dynkin labels $\ell=(\ell_1, \ell_2)$, $m=(m_1, m_2)$ and  $n=(n_1, n_2)$.
After \cite{BerensteinZelevinsky}, the multiplicity of $V_{n}$ in the tensor product $V_\ell \otimes V_m$ is equal to the dimension of the space of invariants
\[
\left(V_{\ell} \otimes V_m \otimes V_n^*\right)^{\SU{3}}.
\]
We set
\begin{equation}\label{def tm}
\tm{\ell}{m}{n} = \dim \left( V_{\ell} \otimes V_m \otimes V_n \right)^{\SU{3}}.
\end{equation}
These integers are called \emph{triple multiplicities for $\SU{3}$}. The multiplicity of  $V_{n}$ in  $V_\ell \otimes V_m$ is thus $\tm{\ell}{m}{n^*}$.

The support of the triple multiplicities (the set of all triples $(\ell, m, n)$ of Dynkin labels such that $\tm{\ell}{m}{n} \neq 0$) generate a sublattice $\lattice{TM}$ of $\ZZ^6$, which is defined by the condition:
\begin{equation}\label{lattice equation}
\ell_1 + m_1 + n_1\equiv \ell_2 + m_2 + n_2 \mod 3.
\end{equation}

\subsection{Linear symmetries}\label{Gg}

A {\em linear symmetry of the triple multiplicities} is a linear automorphism $\theta$ of the lattice $\lattice{\TM}$ that leaves the triple multiplicities unaffected. That is, such that the identity $c(\theta(\ell,m,n)) = \tm{\ell}{m}{n}$ holds.

The definition \eqref{def tm} of the triple multiplicities  makes obvious that the six permutations of the three labels $\ell$, $m$, $n$ are linear symmetries for the triple multiplicities. Another linear symmetry is the ``duality'' symmetry, corresponding to changing each irreducible representation with its dual:
\[
(\ell, m, n) \leftrightarrow (\ell^*, m^*, n^*).
\]
The duality symmetry and the six label permutations generate a group of $12$ linear symmetries isomorphic to $\S{2} \times \S{3}$. We denote this group $\Gg$ and call it the \emph{group of general linear symmetries} for $\SU{3}$, because this group appears also as group of linear symmetries for the triple multiplicities of $\SU{k}$ for any $k$.

\subsection{Polyhedral description of the triple multiplicities}\label{cones}

\begin{figure}[ht] 
	\centering
		\includegraphics[scale=0.3]{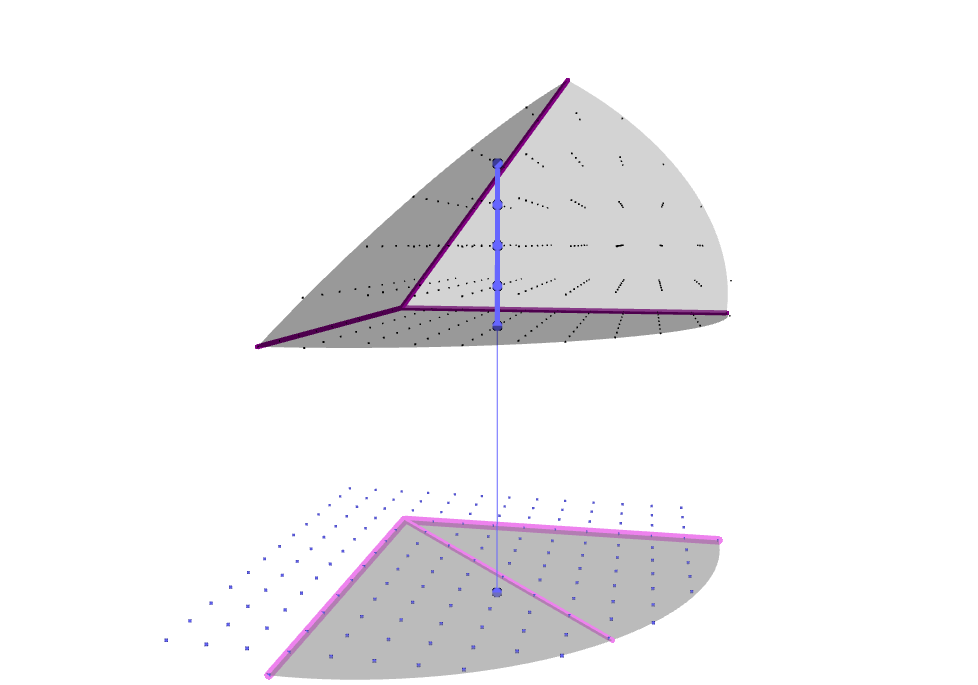}
                \caption{A cone projection and its fiber-counting function. Here a 3-dimensional cone (top) is projected to a 2.dimensional space, with two chambers. The associated fiber-counting function, defined on the lattice points of the 2-dimensional space, counts the lattice points in the corresponding fiber. One such fiber is represented (vertical segment in the top cone). The ``triple multiplicities'' function for $\SU{3}$ is such a fiber-counting function, where the lattice points in the top cone (7-dimensional) are the BZ triangles, and the projection is to a 6-dimensional space (the space of the $(\ell_1, \ell_2; m_1, m_2; n_1, n_2)$).}\label{fig:projection}
\end{figure}

In a finite--dimen\-sion\-al real vector space $\lspace{}$ endowed with a full--rank lattice $\lattice{}$, (a subgroup generated by a basis of $\lspace{}$), a \emph{convex rational polyhedral cone} is the set of all linear combinations, with nonnegative real coefficients, of some fixed finite set of vectors. These cones can also be described as the solution sets of systems of finitely many linear equations $f_i(x) \geq 0$, where the $f_i$ take integer values on the lattice points. In the sequel, \emph{cone} means  \emph{convex rational polyhedral cone}.

A \emph{pointed} cone $C$ in  $\lspace{}$ (cone not containing any line) and a linear projection $\pr{}: \lspace{} \rightarrow \lspace{}'$, sending the lattice $\lattice{}$ onto a lattice $\lattice{}'$ and the cone $C$ to a cone $C'$, defines a \emph{fiber-counting function} on $\lattice{}'$: its value at $t$ is the number of lattice points in $C$ with image $t$. See Figure \ref{fig:projection}.

It follows from \cite{BerensteinZelevinsky} that the function that associates to three irreducible representations of $\SU{3}$ (and, more generally, of $\SU{k}$) the corresponding triple multiplicity (or tensor multiplicity)  is such a fiber-counting function.  In this case, the lattice points in $C$ are combinatorial objects called \emph{Berenstein--Zelevinsky triangles}.

The analogous statement holds for the $\GL{3}$ tensor multiplicities.

A cone $C \subset \lspace{}$ and a projection $\pr{}: \lspace{} \rightarrow \lspace{}'$, mapping $C$ to $C'$, also define a complex of cones subdividing $C'$ 
(a collection of cones, with union $C'$, such that the intersection of any two of them is a face of each; and the faces of any cone in this collection also belong to the collection.). This complex, called the \emph{chamber complex}, has as open cells the sets of points belonging to the projections of exactly the same faces of $C$.
Its maximal closed cells are called the \emph{chamber} of the chamber complex. 
The fiber-counting function associated to $C$ and $\pr{}$  is piecewise quasipolynomial, with the chambers as the domains of validity of the quasipolynomial formulas.

In the cases of the $\SU{k}$ triple multiplicities, or of the $\SU{k}$ or $\GL{k}$ tensor multiplicities, the formulas are actually polynomials (see \cite{Rassart}).

%%%%%%%%%%%%%%%%%%%%%%%%%%%%%%%%%%%%%%%%%%%%%%%%%%%%%%%%%%%%%%%%%%%%%%%%%%%%%%%

\section{Berenstein-Zelevinsky triangles for $\SU{3}$}\label{sec: BZ triangles}

The Littlewood--Richardson coefficients and the $\SU{k}$ triple multiplicities 
can be calculated by means of the Littlewood-Richardson rule and its many avatars,  see \cite{PakVallejo:cones}.
One of these avatars, introduced in \cite{BerensteinZelevinsky}, describes the $\SU{k}$ triple multiplicities as counting combinatorial objects called \emph{Berenstein-Zelevinsky triangles} (``BZ triangles'' in the sequel). These are our main tools in this paper.

In this section, we review the definition of the BZ triangles for the $\SU{3}$ triple multiplicities. (For the general description of the BZ triangles for $\SU{k}$, for any $k$, see  \cite{BerensteinZelevinsky} or \cite{PakVallejo:cones}). These are points in a $7$-dimensional  subspace $\spanBZ$ of a $9$-dimensional space. We calculate a parameterization for the subspace $\spanBZ$ that will be used in the calculations of the next sections.

\subsection{BZ triangles for $\SU{3}$}

\begin{figure}[ht]
  \centering
\begin{tabular}{lr}
 \begin{minipage}[h]{0.55\textwidth}
  \begin{tikzpicture}[scale=0.5, Z1/.style={w, label=below:$\qquad \scriptstyle {Z_1=(2,0,1)}$},
      Z2/.style={w,label=right:$\scriptstyle {Z_2=(2,1,0)}$},
      Z3/.style={w, label=right:$\scriptstyle {Z_3=(1,2,0)}$},
      Z4/.style={w,label=left:$\scriptstyle {Z_4=(0,2,1)}$},
      Z5/.style={w, label=left:$\scriptstyle {Z_5=(0,1,2)}$},
      Z6/.style={w, label=below:$\scriptstyle {Z_6=(1,0,2)}$},
      Y1/.style={w, label=right:$\scriptstyle {Y_1=(3,0,0)}$},
      Y2/.style={w, label=left:$\scriptstyle {Y_2=(0,3,0)}$},
      Y3/.style={w, label=left:$\scriptstyle {Y_3=(0,0,3)}$},
      X/.style=zero]
    \facediagram
  \end{tikzpicture}
  \end{minipage}
  &
  \begin{minipage}{0.37\textwidth}
    \begin{tikzpicture}[scale=0.5]
   \BZT{z_1}{z_2}{z_3}{z_4}{z_5}{z_6}{y_1}{y_2}{y_3}
    \draw (Z6) -- (Y3) node [midway, scale=0.7, below] () {$m_1$};
    \draw (Z1) -- (Y1) node [midway, scale=0.7, below] () {$m_2$};
    \draw (Z2) -- (Y1) node [midway, scale=0.7, right] () {$n_1$};
    \draw (Z3) -- (Y2) node [midway, scale=0.7, right] () {$n_2$};
    \draw (Z4) -- (Y2) node [midway, scale=0.7, left] () {$\ell_1$};
    \draw (Z5) -- (Y3) node [midway, scale=0.7, left] () {$\ell_2$};
    \end{tikzpicture}
    
    \vspace{1mm}
    \end{minipage}
\end{tabular}
    \caption{The BZ Graph $\Gamma$ (left) and the coordinates of the BZ labellings (right).}\label{BZ graph}%\label{BZtriangle}
\end{figure}
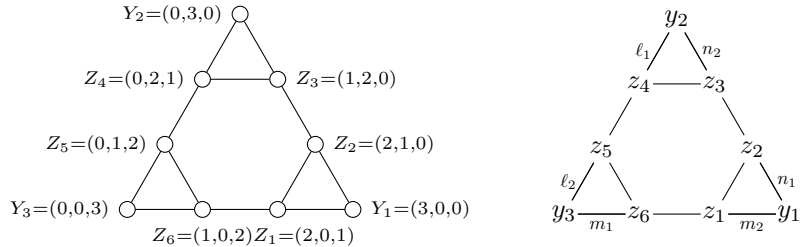

Consider the \emph{BZ graph} $\Gamma$, shown in Figure \ref{BZ graph}: its vertices are the 9 points $(i,j,k)$ with nonnegative integer coordinates fulfilling $i+j+k=3$, different from $(1,1,1)$. There is an edge between any two vertices with difference  $(1,-1,0)$, $(1,0,-1)$ or $(0,1,-1)$. The 9 vertices of the BZ graph are the 3 vertices of an equilateral triangle and the 6 vertices of a regular hexagon inscribed in the triangle. We will refer to the vertices of the BZ graph as $Y_1, Y_2, Y_3$ and $Z_1, \ldots, Z_6$ as shown in Figure \ref{BZ graph}.

For any labelling of the BZ graph, we will denote $y_1, y_2, y_3$ and $z_1, \ldots, z_6$ the labels of the vertices  $Y_1, Y_2, Y_3$ and $Z_1, \ldots, Z_6$ (see Figure \ref{BZ graph}). In the $9$--dimensional space $\RR^{\Gamma}$ of all BZ graph labellings, let $\spanBZ$ be the $7$--dimensional subspace defined by the equations: 
\begin{equation}\label{BZ condition}
  z_1-z_4=z_5-z_2=z_3-z_6.
\end{equation}
The points of $\spanBZ$ are the BZ graph labellings such that any side of the hexagon sums as much as the opposite side. Let $\coneBZ$ be the cone of all  points in $\spanBZ$ with \emph{nonnegative} labels. Let $\latBZ$ be the set of all integer points in the cone $\BZ$. The \emph{Berenstein--Zelevinsky triangle} (\emph{BZ triangle} in the sequel) are the elements of $\latBZ$; otherwise said, the BZ triangles are the labellings of the BZ graphs with nonnegative integer labels, fulfilling \eqref{BZ condition}.

Let us introduce also the linear map $\proj: \spanBZ \rightarrow \RR^6$ that sends a BZ graph labelling in $\spanBZ$ to the tuple $(\ell_1, \ell_2; m_1, m_2; n_1, n_2)$ defined by 
\begin{equation}\label{coordinates from BZ}
  \begin{array}{lll}
    \ell_1 = y_2 + z_4, & m_1 = y_3 + z_6, & n_1 = y_1 + z_2,\\
    \ell_2 = y_3 + z_5, & m_2 = y_1 + z_1, & n_2 = y_2 + z_3.
    \end{array}
\end{equation}
These are the sums of the labels at some vertex $Y_i$ and at one of the two neighbouring vertices (see again Figure \ref{BZ graph}). The linear map $\proj$ maps the integer points in $\spanBZ$ onto the lattice $\lattice{TM} \subset \ZZ^6$ of all points $(\ell_1, \ell_2, m_1,m_2, n_1, n_2)$ fulfilling \eqref{lattice equation}.
%\[
%\ell_1 + m_1 + n_1 \equiv \ell_2 + m_2 + n_2 \mod 3.
%\]
Then the triple multiplicity for the weights $\ell=(\ell_1, \ell_2)$, $m=(m_1, m_2)$ and $n = (n_1, n_2)$ counts the BZ triangles in the fiber over $t = (\ell_1, \ell_2; m_1, m_2; n_1, n_2)$ of the projection $\proj$:
\[
\tm{\ell}{m}{n} = \# \left(\proj^{-1}(\ell;m;n) \cap \latBZ\right).
\]

\begin{table}
  \centering
    \caption{The objects in polyhedral geometric description of the $\SU{3}$ Triple Multiplicities.}

  \begin{tabular}{lp{10cm}}
    \toprule
    $\lspace{TM}$ & 6D vector space of the $(\ell_1, \ell_2; m_1, m_2; n_1, n_2)$. \\
    $\lattice{TM}$ & Lattice of the integer points in $\lspace{TM}$ fulfilling $\ell_1+m_1 + n_1 \equiv \ell_2+m_2+n_2 \mod 3$. \\
    $\TM$ & Cone in $\lspace{TM}$. \\
    $ \latTM$ &  Lattice points in $\TM$. Support of the triple multiplicities. \\
    $\chK$ & Chamber complex for the triple multiplicities.  Subdivides $\TM$.  \\[2mm]
    $\spanBZ$ & 7D vector space of all BZ Graph labellings fulfilling \eqref{BZ condition}. Endowed with a projection onto $\lspace{TM}$.\\
    $\lattice{BZ}$ & Lattice of all integer labellings in $\spanBZ$. Projects onto $\lattice{TM}$.\\
    $\BZ$ & Cone of all nonnegative labellings in $\spanBZ$. Projects onto $\TM$.\\
    $\latBZ$ &  Set of all  Berenstein-Zelevinsky Triangles. (All lattice points in $\BZ$). \\
    \bottomrule
  \end{tabular}
\end{table}

\subsection{A parameterization of the BZ triangles}\label{parameterization}

We now parameterize the space $\spanBZ$. We  use as parameters the values of $\ell_1$, $\ell_2$, $m_1$, $m_2$, $n_1$, $n_2$ and the label $x=-y_1$ of vertex $Y_1$.

From \eqref{coordinates from BZ} we get
\begin{equation}\label{z equal}
    \begin{array}{lll}
    z_4 = \ell_1 - y_2, & z_6 = m_1 - y_3, & z_2 = n_1 - y_1,\\
    z_5 = \ell_2 - y_3, & z_1 = m_2 - y_1, & z_3 = n_2 - y_2.
    \end{array}
\end{equation}

We set $\omega=z_4-z_1$. After \eqref{BZ condition}, there is also $\omega=z_6-z_3=z_2-z_5$. Averaging these three expressions for $\omega$ gives
\begin{equation}\label{omega z}
\omega = \frac{1}{3}(z_2+z_4+z_6-z_1-z_3-z_5).
\end{equation}
Replacing in  \eqref{omega z}
each $z_i$ with its expression in terms of $t=(\ell_1, \ell_2; m_1,m_2;n_1,n_2)$ and the $y_j$ from \eqref{z equal}
yields
\[
\omega = \frac{1}{3}(\ell_1+m_1+n_1-\ell_2-m_2-n_2)
\]
(the $y_i$ cancel out).

Replacing in \eqref{BZ condition} each $z_i$ with its expression in \eqref{z equal} 
yields 
three relations $y_i = y_1+f_i(t)$, with the linear forms $f_i(t)$ displayed in Table \ref{table:fi and gj}, where
\[
  \omega(t)=\frac{1}{3}\left(\ell_1+m_1+n_1-\ell_2-m_2-n_2\right).
  \]

Replacing $y_i$ with $y_1+f_i(t)$ in \eqref{z equal} yields six relations $z_j = -y_1-g_j(t)$, with the forms $g_j(t)$ also shown 
in Table \ref{table:fi and gj}.

\begin{figure}[ht]
  \centering
  \begin{tikzpicture}[scale=0.6]
   \input{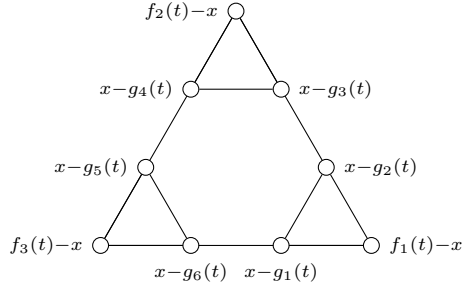}
    \end{tikzpicture}
  \caption{The labelled BZ graph $BZ(t,x)$. The linear forms $f_i(t)$ and $g_j(t)$ are those defined in Table \ref{table:fi and gj}.}\label{BZtx}
\end{figure}

We finally set $x=-y_1$. For any $t=(\ell_1, \ell_2; m_1, m_2; n_1, n_2)$ and $x$, let $BZ(t,x)$ be the labelled BZ graph shown in Figure \ref{BZtx}.
Then the linear map $(t,x) \mapsto BZ(t,x)$ establishes an isomorphism from $\RR^6 \times \RR$ to the linear span $\spanBZ$ of the BZ triangles.
This isomorphism of real vector spaces restricts to an isomorphism of lattices from $\Lambda \times \ZZ$ to the lattice $\lattice{BZ}$ of the integer points of $\spanBZ$.

Since $y_i=y_1+f_i(t)$ and $x=-y_1$, the inequality $y_i \ge 0$ is equivalent to $f_i(t) \ge x$.
Similarly, since $z_j=-y_1-g_j(t)$, the inequality $z_j \ge 0$ is equivalent to $g_j(t) \le x$.
Therefore, under the parameterization considered here, the cone $\coneBZ$ is the set of solutions of
\[
\begin{cases}
  \forall i, x \le f_i(t),\\
    \forall j, x \ge g_j(t).
\end{cases}
\]

\begin{table}
  \centering
  \caption{The linear forms $f_i$ and $g_i$.}\label{table:fi and gj}
 \[
 \begin{array}{l@{\quad}l@{\quad}l}
  f_{1}(t)  = 0,
  & f_{2}(t)  = \ell_1-m_2-\omega(t),
  &f_3(t) = \ell_2-n_1+\omega(t),\\
  g_1(t) = -m_2,
  & g_3(t) = \ell_1-m_2-n_2-\omega(t),  
  & g_5(t) = -n_1+\omega(t),\\
  g_2(t) = -n_1,
  &g_4(t) = -m_2-\omega(t),
  &g_6(t) = \ell_2-m_1-n_1+\omega(t)  \\[3mm]
  \multicolumn{3}{c}{\text{with } \omega(t)=\frac{1}{3}\left(\ell_1+m_1+n_1-\ell_2-m_2-n_2\right).}
\end{array}
\]
\end{table}

%%%%%%%%%%%%%%%%%%%%%%%%%%%%%%%%%%%%%%%%%%%%%%%%%%%%%%%%%%%%%%%%%%%%%%%%%%%%%%%

\section{Linear symmetries}\label{sec:symmetries}

In this section, we determine all linear symmetries for the triple multiplicities.
As a first step, in \ref{subsec:symBZ}, we find the group $\upgroup$ of all linear symmetries for the set $\latBZ$ of all BZ triangles.
These are all automorphisms of $\lattice{BZ}$ that stabilize the set $\latBZ$.
Then, in \ref{sec:liftable},  we check that each element $\theta$ of  $\upgroup$ induces a linear symmetry $\theta'$ of the triple multiplicities.
In such a situation, we say that \emph{$\theta'$ lifts to $\theta$}, or that  \emph{$\theta'$ admits $\theta$ as a lifting.}
We obtain this way a group of symmetries $\Gl$ of the triple multiplicities, isomorphic to $\upgroup$, of order $72$: the group of symmetries ``that admit a lifting''.
But not all symmetries of the triple multiplicities admit a lifting: the duality symmetry does not.
The group of all linear symmetries, $G$, is thus bigger than $\Gl$. 
In \ref{sec:all symmetries}, we embed $G$ into a permutation group (the group of permutations of all rays of the cone $\TM$) that is found to have order $144$. This is enough to conclude that $G$ has exactly $144$ elements, and is generated by $\Gl$ and the duality symmetry. In the course of the proof, we determine all rays of the cones $\BZ$ and $\TM$.

As subgroup of interest, we also consider, besides $\Gg$ (group of general symmetries, see \ref{Gg}) the intersection  $\Glg=\Gg \cap \Gl$. Figre \ref{inclusions} shows the relations between these groups (see \ref{sec:liftable} for the fact that the map $\upgroup \rightarrow \Gl$ is an ismorphism).

      \tikzset{
  symbol/.style={
    draw=none,
    every to/.append style={
      edge node={node [sloped, allow upside down, auto=false]{$#1$}}}
  }
      }
      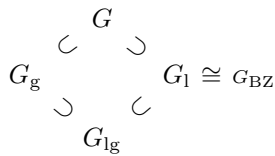
\begin{figure}[ht]
      \centering
      \begin{tikzcd}[column sep=.7em, row sep =.7 em]
            &  G    &                 &\\
        \Gg \arrow[ur,symbol=\subset]&       & \Gl \arrow[ul ,symbol=\subset] \arrow[r, symbol=\cong]   & \scriptstyle{\upgroup} \\
            &  \Glg \arrow[ul,symbol=\subset] \arrow[ur,symbol=\subset] &  &\\
      \end{tikzcd}
      \caption{Relations between the groups of symmetries.}\label{inclusions}
\end{figure}

The embedding of $G$ into a group of permutations allows to elucidate the structure of the groups $\Gg$, $\Gl$ and $\Glg$. The results are detailed in \ref{subsec:sumsym} and summarized in Table \ref{fig:subgroups}.

      \begin{table}[ht]
        \centering
             \caption{Groups of linear symmetries for the triple multiplicities.}\label{fig:subgroups}
     \begin{tabular}{lrc}
    \toprule
    Group & Order   & Structure  \\
    \midrule
    $G$ (all linear symmetries)  & 144 & $\S{2} \times \left(\S{3} \wr \S{2}\right)$ \\
    $\Gg$ (general symmetries)   &  12 & $\S{3} \times \S{2}$ \\
    $\Gl$ (symmetries that lift) &  72 & $\S{2} \times \S{3} \times \S{3}$ \\
    $\Glg$ ($=\Gl \cap \Gg$)     &   6 & $\S{3}$     \\
    \bottomrule
  \end{tabular}  
     \end{table}

\subsection{Symmetries for the BZ triangles of $SU(3)$}\label{subsec:symBZ}

\subsubsection{General symmetries for the BZ triangles}\label{subsubsec:genB}

\begin{figure}[ht]
  \centering
  \begin{tikzpicture}[scale=0.6, Z1/.style={w, label=below:${z_1}$},
      Z2/.style={w,label=right:${z_2}$},
      Z3/.style={w, label=right:${z_3}$},
      Z4/.style={w,label=left:${z_4}$},
      Z5/.style={w, label=left:${z_5}$},
      Z6/.style={w, label=below:${z_6}$},
      Y1/.style={w, label=below:${y_1}$},
      Y2/.style={w, label=left:${y_2}$},
      Y3/.style={w, label=below:${y_3}$},
      X/.style=zero]
    \facediagram
        \draw (Z6) -- (Y3) node [midway, scale=0.7, below] () {$m_1$};
    \draw (Z1) -- (Y1) node [midway, scale=0.7, below] () {$m_2$};
    \draw (Z2) -- (Y1) node [midway, scale=0.7, right] () {$n_1$};
    \draw (Z3) -- (Y2) node [midway, scale=0.7, right] () {$n_2$};
    \draw (Z4) -- (Y2) node [midway, scale=0.7, left] () {$\ell_1$};
    \draw (Z5) -- (Y3) node [midway, scale=0.7, left] () {$\ell_2$};
         \node (sw) at (210:1) {};  
         \path (Y1) --++(-30:1) node  (se) {};
         \draw (sw) -- ++(30:8cm) [color=red, dash pattern=on 2pt off 3pt on 4pt off 4pt] node (ne) {$s_1$};
         \draw (se) -- ++(150:8cm) [color=blue, dash pattern=on 2pt off 3pt on 4pt off 4pt] node (nw) {$s_2$};
         \path (ne) --++(210:0.3) node (se arrow) [rotate=120, color=red] {$\longleftrightarrow$};
         \path (nw) --++(-30:0.3) node (nw arrow) [rotate=60, color=blue] {$\longleftrightarrow$};
       \end{tikzpicture}  
  \caption{Action of the simple transpositions $s_1$ and $s_2$ on the BZ triangles and on the triple multiplicities.}\label{BZtriangle symmetries}
\end{figure}
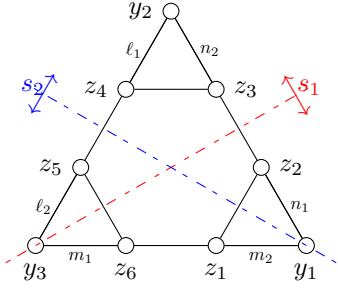

The BZ triangles have been introduced in \cite{BerensteinZelevinsky} to make obvious some of the general symmetries of the $\SU{3}$ tensor multiplicities. Namely, the 6 linear symmetries of the triangle $Y_1 Y_2 Y_3$ preserve the whole BZ graph, and permute the coordinates $y_i$ and $z_j$ of the BZ graph labellings, in such a way that the relations \eqref{BZ condition} are preserved. These six symmetries induce therefore symmetries of $\latBZ$. They also induce symmetries of the triple multiplicities (since these are obtained as sums of the labels $y_i$ and some neighboring $z_j$) forming a group $\Glg$. The group of the six linear symmetries of the triangle is generated by the reflections $s_1$ and $s_2$ with respect to the bisectors through $Y_1$ and $Y_3$. The effect on these generators on the BZ triangles and on the triple multiplicities is easily read on Figure \ref{BZtriangle symmetries}.  The groups of the symmetries of the triangle, of the induced symmetries of the BZ triangles, and $\Glg$ (induced symmetries of the triple multiplicities) are all isomorphic, and isomorphic to the symmetric group $\S{3}$.

Remember that the 12 general linear symmetries of the triple multiplicities for $SU(3)$ form a group $\Gg$ isomorphic to $\S{2} \times \S{3}$, where the factor $\S{2}$ is generated by the duality involution $(\ell; m; n) \leftrightarrow (\ell^*; m^*; n^*)$, and the elements of the factor $\S{3}$ are the permutations of the triple $(\ell; m; n)$.
The elements of $\Glg$  are the 3 even permutations of $(\ell;m;n)$, together with the 3 transpositions composed with the duality involution. Otherwise stated, $\Glg$ is generated by $(\ell;m;n) \leftrightarrow (m^*; \ell^*; n^*)$ and by  $(\ell;m;n) \leftrightarrow (\ell^*; n^*; m^*)$. This subgroup $\Glg$ is isomorphic to $\S{3}$, but distinct from the group of the permutations of $(\ell;m;n)$.

Note that all of the above holds not only for the $SU(3)$ case considered here, but also for $SU(k)$ for any $k$ (see \cite[Remark (a) p.10]{BerensteinZelevinsky}).

\subsubsection{All symmetries of the BZ triangles of $SU(3)$}\label{subsubsec:symBZ}

As already mentioned, the BZ triangles generate a cone $\coneBZ$ in the subspace $\spanBZ$ of $\RR^{\Gamma}$ defined by the equations $z_1-z_4 =z_3-z_6=z_5-z_2$. We will determine its rays.

For any $T \in \Lambda$ whose triple multiplicity is $1$, i.e. which is the projection of a unique BZ triangle. Denote this BZ triangle by $\Delta_T$.  We contend that 
the rays of $\coneBZ$  are generated by the eight BZ triangles  $\upray{C_1}$, $\upray{C_2}$, $\upray{C_3}$, $\upray{D_3}$, $\upray{D_5}$, $\upray{D_1}$, $\upray{\ltriangle}$ and $\upray{\rtriangle}$ shown in Table \ref{fundamental BZT}. Following \cite{CoquereauxZuber}, we call them the \emph{fundamental BZ triangles}.

\begin{table}[ht]
  \centering
    \caption{The eight fundamental BZ triangles $\upray{C_i}$, $\upray{D_i}$, $\upray{\ltriangle}$ and $\upray{\rtriangle}$ of the cone $\coneBZ$, with their projections $\ray{C_i}$, $\ray{D_i}$, $\ray{\ltriangle}$ and $\ray{\rtriangle}$.}\label{fundamental BZT}
  \[
  \begin{array}{c@{\qquad}c@{\qquad}c}
    \toprule
    \begin{tikzpicture}[scale=0.35] \fullBZT{0}{0}{0}{0}{0}{0}{1}{0}{0}\end{tikzpicture}
    &
    \begin{tikzpicture}[scale=0.35] \fullBZT{0}{0}{0}{0}{0}{0}{0}{1}{0}\end{tikzpicture}
    &
    \begin{tikzpicture}[scale=0.35] \fullBZT{0}{0}{0}{0}{0}{0}{0}{0}{1}\end{tikzpicture}
    \\
    \upray{C_1}& \upray{C_2} & \upray{C_3}\\
    \ray{C_1}=(00|01|10) &  \ray{C_2}=(10|00|01)  &  \ray{C_3}=(01|10|00) \\[5mm]
     \begin{tikzpicture}[scale=0.35] \fullBZT{0}{0}{1}{0}{0}{1}{0}{0}{0}\end{tikzpicture}
    &
    \begin{tikzpicture}[scale=0.35] \fullBZT{0}{1}{0}{0}{1}{0}{0}{0}{0}\end{tikzpicture}
    &
    \begin{tikzpicture}[scale=0.35] \fullBZT{1}{0}{0}{1}{0}{0}{0}{0}{0}\end{tikzpicture}
    \\
    \upray{D_3} &  \upray{D_5} &  \upray{D_1} \\
    \ray{D_3}=(00|10|01) & \ray{D_5}=(01|00|10) & \ray{D_1}=(10|01|00)\\[5mm]
    \begin{tikzpicture}[scale=0.35] \fullBZT{1}{0}{1}{0}{1}{0}{0}{0}{0}\end{tikzpicture}
    & &
    \begin{tikzpicture}[scale=0.35] \fullBZT{0}{1}{0}{1}{0}{1}{0}{0}{0}\end{tikzpicture}
    
    \\
    \upray{\ltriangle} & &\upray{\rtriangle} \\
    \ray{\ltriangle}=(01|01|01) & &\ray{\rtriangle}=(10|10|10)\\
    \bottomrule
  \end{array}
  \]
  \end{table}

Indeed, given a point $T$ in $\coneBZ$, with labels $y_i$ and $z_j$ (as in Figure \ref{BZ graph}), set (as in \ref{parameterization}) $\omega~=~z_4~-~z_1~=z_6-z_3=z_2-z_5$. If $\omega \ge 0$, then $T$ decomposes as 
\begin{equation}\label{upray coordinates}
  T= y_1 \cdot \upray{C_1} + y_2 \cdot \upray{C_2} +  y_3 \cdot \upray{C_3}
  + z_1 \cdot \upray{D_1} + z_3 \cdot \upray{D_3} + z_5 \cdot \upray{D_5} +
  \omega \cdot \upray{\rtriangle}
\end{equation}
with nonnegative coefficients; and if $\omega \le 0$, 
\begin{equation}\label{alt upray coordinates}
  T= y_1 \cdot \upray{C_1} + y_2 \cdot \upray{C_2} +  y_3 \cdot \upray{C_3}
  + z_4 \cdot \upray{D_1} + z_6 \cdot \upray{D_3} + z_2 \cdot \upray{D_5} +
  (-\omega) \cdot \upray{\ltriangle}.
\end{equation}
This shows that the fundamental BZ triangles generate $\coneBZ$. These are $8$ vectors spanning the $7$-dimensional space $\spanBZ$. There is therefore a single linear relation between them, which is
\begin{equation}\label{rel between fundamental}
      \upray{D_1} + \upray{D_3} + \upray{D_5}=
  \upray{\ltriangle} + \upray{\rtriangle}.
\end{equation}
Therefore, none of the fundamental BZ triangles is a positive linear combination of the other. This completes the checking of the fact that the fundamental triangles are generators of the rays of $\coneBZ$.

Note that the fundamental BZ triangles do not generate only the vector space~$\spanBZ$, but the lattice of the integer points of $\spanBZ$, since any point $T$ in this lattice admits a decomposition \eqref{upray coordinates} with integer coefficients.

The linear symmetries of $\latBZ$ (the linear automorphisms of $\lattice{BZ}$ that stabilize $\latBZ$) form a group $\upgroup$.
Any element $\theta$ of $\upgroup$ permutes the fundamental BZ triangles, since they are the minimal ray generators of the cone $\coneBZ$.
Taking into account the relation \eqref{rel between fundamental} between the fundamental BZ triangles, $\theta$ stabilizes the set $\{\upray{D_3}, \upray{D_5},\upray{D_1}\}$ (the only set of 3 fundamental BZ triangles whose sum is the sum of two other fundamental BZ triangles).
Similarly, $\theta$ stabilizes $\{\upray{\ltriangle},\upray{\rtriangle}\}$ (the only set of two fundamental BZ triangles whose sum is the sum of three other fundamental BZ triangles). Finally, $\theta $ also stabilizes the set of the 3 remaining rays $\{\upray{C_1}, \upray{C_2},\upray{C_3}\}$.
Reciprocally, any permutation of the 8 fundamental BZ triangles that stabilizes  each of $\{\upray{C_1}, \upray{C_2},\upray{C_3}\}$, $\{\upray{D_3}, \upray{D_5},\upray{D_1}\}$ and  $\{\upray{\ltriangle},\upray{\rtriangle}\}$ lifts to a unique element of $\upgroup$, since it leaves the relation \eqref{rel between fundamental} unchanged.
This shows that $\upgroup$ is isomorphic to
\[
\S{\{\upray{C_1}, \upray{C_2},\upray{C_3}\}} \times \S{\{\upray{D_3}, \upray{D_5},\upray{D_1}\}} \times  \S{\{\upray{\ltriangle},\upray{\rtriangle}\}}.
\]
In particular, $\upgroup$ has order $3!\times3!\times2!=72$.

The elements of $\upgroup$ can be described as linear symmetries of parts of the BZ graphs. The permutations of $\{\upray{C_1}, \upray{C_2},\upray{C_3}\}$ are obtained by applying the linear symmetries of the triangle $Y_1 Y_2 Y_3$, only to $Y_1$, $Y_2$, $Y_3$ (leaving unaffected the vertices~$Z_j$). The elements of
$\S{\{\upray{D_3}, \upray{D_5},\upray{D_1}\}} \times  \S{\{\upray{\ltriangle},\upray{\rtriangle}\}}$ are the linear symmetries of the hexagon $Z_1 Z_2 Z_3 Z_4 Z_5 Z_6$  (dihedral group of order $12$) applied to its vertices (leaving unaffected the vertices $Y_i$).

Alternatively, $\upgroup$ is generated by the $6$ linear symmetries of the BZ graph (the linear symmetries of the triangle $Y_1 Y_2 Y_3$, applied to all $9$ vertices of the BZ graph) and the $12$ linear symmetries of the hexagon $Z_1 Z_2 Z_3 Z_4 Z_5 Z_6$ applied only to the vertices of this hexagon.

\subsection{The linear symmetries of the triple multiplicities that lift to linear symmetries of the set of all BZ triangles}\label{sec:liftable}

Any element $\theta$ of $\upgroup$ is, a priori, of the form $BZ(t,x) \mapsto BZ(t',x')$. The symmetry $\theta \in \upgroup$ induces a linear symmetry of the triple multiplicities if and only if $t'$ is independent of $x$. Let us show that this is always the case. Let $y_i$ and $z_j$ be the coordinates of  $BZ(t,x)$, and let $y'_i$ and $z'_j$ be the coordinates of $BZ(t',x')$. Then $\theta$
permutes the coordinates $y_i$, and also permutes the coordinates $z_j$. After the relations \eqref{coordinates from BZ}, the coordinates of $t'$ are all of the form $y'_i+z'_j$. Each of them is thus $y_p+z_q$ for some $p$ and $q$.  But $y_p+z_q=(f_p(t)-x)+(x-g_q(t))$, which is independent on $x$. This proves that $t'$ is independent on $x$, as announced.
As a conclusion, any element $\theta$ of $\upgroup$ induces a linear symmetry of the triple multiplicities.

A direct calculation shows that the 8 fundamental BZ triangles have 8 distinct projections  $\ray{C_i}$, $\ray{D_i}$, $\ray{\ltriangle}$ and $\ray{\rtriangle}$ (given in Table \ref{fundamental BZT}). As a consequence, the
symmetries of the triple multiplicities induced by the 72 elements of $\upgroup$ are all distinct. They form a subgroup $\Gl$ of the group of all symmetries of the triple multiplicities. The subgroup $\Gl$ is isomorphic to $\upgroup$, and thus to $\S{2} \times \S{3} \times \S{3}$, and has order $72$. It can be characterized as the group of all linear automorphisms of the lattice $\Lambda$ that stabilizes each of the sets $\{\ray{C_1},\ray{C_2},\ray{C_3}\}$,  $\{\ray{D_1},\ray{D_3},\ray{D_5}\}$ and  $\{\ray{\ltriangle}, \ray{\rtriangle}\}$.  

The group $\Gl$ contains the group $\Glg$ of the six ``general'' symmetries that admit a lifting, which are easily recognized by their effect on the projections  $\ray{C_i}$, $\ray{D_i}$, $\ray{\ltriangle}$ and $\ray{\rtriangle}$. Namely, $\Glg$ is the subgroup of the elements that permute the pairs $(\ray{C_1}, \ray{D_5})$, $(\ray{C_2},\ray{D_3})$, $(\ray{C_3},\ray{D_1})$, and, only for the elements of order $2$, swap  $\ray{\ltriangle}$ with $\ray{\rtriangle}$.

The duality involution $(\ell_1, \ell_2; m_1,m_2; n_1,n_2) \leftrightarrow  (\ell_2, \ell_1; m_2,m_1; n_2,n_1)$, and the odd permutations of $(\ell; m; n)$ do not lift, since they do not stabilize $\{\ray{C_1},\ray{C_2},\ray{C_3}\}$. See the final remark \ref{no lifting} for a more straightforward argument for this fact.

\subsection{All linear symmetries of the $\SU{3}$ triple multiplicities}\label{sec:all symmetries}

Remember that $G$ is the group of all linear symmetries of the triple multiplicities. Any element of $G$ stabilizes the cone $\TM$; therefore it permutes the rays of $\TM$, and also their minimal ray generators. These minimal ray generators are among the elements of
\[
R=\{\ray{C_1}, \ray{C_2}, \ray{C_3}, \ray{D_1}, \ray{D_3}, \ray{D_5}, \ray{\ltriangle}, \ray{\rtriangle}\}.
\]
These vectors span not only the vector space $\lspace{\TM}=\RR^6$, but also the lattice $\lattice{\TM}$, since they are the projections of the fundamental BZ triangles that span the lattice $\lattice{BZ}$ of the integer points in $\spanBZ$.

The following relations between the vectors in $R$ hold:
\begin{equation}\label{rel between R}
\ray{C_1}+\ray{C_2}+\ray{C_3}
=\ray{D_1}+\ray{D_3}+\ray{D_5}
=\ray{\ltriangle}+\ray{\rtriangle}.
\end{equation}
Both can be checked straightforwardly. (Actually the second relation is obtained from \eqref{rel between fundamental} by projecting. For the fact that the other relation is \emph{not} obtained by such a projection, see Section \ref{no lifting}). Since $R$ is a set of 8 vectors spanning $\RR^6$, there exist two independent relations between the elements of $R$, which are just \eqref{rel between R}. In particular, there cannot be any relation expressing one vector from $R$ as a positive linear relation of others. 
This shows that $R$ is exactly the set of all minimal ray generators of $\TM$.

We can thus embed $G$ into the group of permutations of $R$. Moreover, any element of $G$ stabilizes the subset $\{\ray{\ltriangle}, \ray{\rtriangle}\}$ (the only pair of vectors of $R$ whose sum is sum of three other vectors of $R$). Similarly, any element of $G$ stabilizes each of $\{\ray{C_1}, \ray{C_2},\ray{C_3}\}$ and  $\{\ray{D_1}, \ray{D_3},\ray{D_5}\}$, or swaps them. Indeed, these two sets are the only sets of three vectors of $R$ whose sums are also sums of two vectors of $R$.

As a conclusion, the embedding of $G$ into the group of permutations of $R$ takes its value in the subgroup $H$ of the group of all permutations of $R$,  that stabilizes $\{\ray{\ltriangle}, \ray{\rtriangle}\}$ and $\{\{\ray{C_1}, \ray{C_2},\ray{C_3}\},\{\ray{D_1}, \ray{D_3},\ray{D_5}\}\}$. This subgroup is isomorphic to $\S{2} \times \left( \S{3} \wr \S{2}\right)$, where $\wr$ denotes the wreath product of groups  (see for example \cite[page 187]{DummitFoote}).  The subgroup $H$  thus has order $2 \times (2 \times (3!)^2)$, which is $144$. On the other hand, $G$ contains the subgroup $\Gl$ of order $72$ consisting of the symmetries that admit a lifting. Therefore, $G$ has order some multiple of $72$. The group $G$ also contains some elements not in $\Gl$, such as the duality automorphism.  Therefore $G$ has order at least $144$. This proves that $G$ has order $144$ exactly, and is isomorphic to $H$.

\subsection{The subgroups of symmetries as permutation groups}\label{subsec:sumsym}

Let us summarize and complete the results in this section. Firstly, the group $G$ of all linear symmetries of the triple multiplicities of~$SU(3)$ has order $144$ and is isomorphic to $\S{2} \times \left(\S{3} \wr \S{2}\right)$. It embeds into the group $\S{R}$ of the permutations of $R$ as the stabilizer of $\{\{\ray{C_1}, \ray{C_2},\ray{C_3}\}, \{\ray{D_1}, \ray{D_3},\ray{D_5}\}\}$.

The subgroup $\Gl$ of all symmetries that lift to linear symmetries of $\latBZ$ has order $72$ and embeds into $\S{R}$ as
$$\S{\{\ray{C_1},\ray{C_2},\ray{C_3}\}} \times \S{\{\ray{D_1},\ray{D_3},\ray{D_5}\}} \times \S{\{\ray{\ltriangle}, \ray{\rtriangle}\}}.$$
Each element of $\Gl$ lifts uniquely to a linear symmetry of $\latBZ$, and therefore $\Gl \cong \upgroup$.

Let us consider the group $\Gg$ of the 12 ``general symmetries'',generated by the permutations of $(\ell;m;n)$ and the duality involution $(\ell;m;n) \leftrightarrow (\ell^*; m^*; n^*)$.
The permutations of $(\ell;m;n)$ permute the pairs $(C_1,D_3)$, $(C_2, D_5)$ and $(C_3, D_1)$, while the duality involution swaps the terms of each of these pairs (i.e. swaps $(C_1, C_2, C_3)$ with $(D_3, D_5, D_1)$), and additionally swaps $\ray{\rtriangle}$ with $\ray{\ltriangle}$.

Finally, the subgroup $\Glg=\Gl \cap \Gg$, made of the 6 ``general symmetries with a lifting'', is generated by the two involutions $s_1$ swapping $(\ell;m;n)$ with $(m^*;\ell^*;n^*)$, and $s_2$ swapping $(\ell;m;n)$ with $(\ell^*;n^*;m^*)$. The involution $s_1$ swaps $\ray{D_3}$ with $\ray{D_5}$, $\ray{C_1}$ with $\ray{C_2}$ and $\ray{\rtriangle}$ with $\ray{\ltriangle}$.
The other generator $s_2$ of $\Glg$ swaps $\ray{D_5}$ with $\ray{D_1}$, $\ray{C_2}$ with $\ray{C_3}$ and $\ray{\rtriangle}$ with $\ray{\ltriangle}$.

%% SU(3)
%%----------------------------------------------------------------------

\section{The chamber complex for the triple multiplicities for \(\SU{3}\)}\label{sec:chamber complex}

In this section. we show that the support of the triple multiplicities is the set of all lattice points of a cone $\TM$, covered by $18$ smaller cones $C(i, j)$, that are the domains for linear formulas for the triple multiplicities. The cones $C(i, j)$  generate a chamber complex $\chK$ subdividing $\TM$. We check that the $C(i, j)$ are full-dimensional, and thus that they are the chambers of $\chK$. We also show that all $C(i,j)$ are simplicial (i.e. have exactly as many facets as the dimension of the ambient space, which is $6$). Then we encode combinatorially all cells of the chamber complex.
As an aside, we are able to count the cells of each dimension.
Finally, we associate to each cell a diagram that makes clear the action of the group $G$ of linear symmetries of the triple multiplicities, and its subgroups, on the cells.

\subsection{The chamber complex $\chK$}

We have seen at the end of Section \ref{parameterization} that the cone $\coneBZ$ is defined by a system of $18$ inequalities:
\[
\begin{cases}
  \forall i \in\{1,2,3\},        & x \le f_i(t),\\
    \forall j \in\{1,2,3,4,5,6\},& x \ge g_j(t),
\end{cases}
\]
with the $3$ forms $f_i$ and the $6$ forms $g_j$ defined in Table \ref{table:fi and gj}. This system can be summarized by the condition:
\begin{equation}\label{BZ polytope}
 \max_q g_q(t) \le x \le \min_p f_p(t).
 \end{equation}
 As a consequence, the labels $t=(\ell;m;n)$ of the non-zero triple multiplicities all belong to the subset $\TM$ of $\lspace{TM}$  defined by
 \[
 \max_q g_q(t) \le \min_p f_p(t).
 \]
 The above inequality summarizes the system of $18$ linear inequalities:
 \[
 \forall i \in \{1,2,3\}, \; \forall j \in \{1,2,3,4,5,6\}, \; g_j(t) \le f_i(t).
 \]
 This is why $\TM$ is a cone.

For any fixed $t=(\ell;m;n)\in \latof{\TM}$,  the triple multiplicity $c(t)$ counts the integers $x$ that fulfill \eqref{BZ polytope}. Therefore
\[
c(t) = 1+ \max(0, \min_p f_p(t) - \max_q g_q(t)).
\]
This shows that the triple multiplicity function is piecewise polynomial of degree $1$ (i.e. piecewise linear with constant term). For each $i\in \{1,2,3\}$ and each $j \in \{1,2,\ldots,6\}$, the linear formula
\[
c(t) = 1+ f_i(t) - g_j(t)
\]
holds for all lattices points $t$ in the set $C(i,j)$ defined by 
  \[
  \min_p f_p(t) = f_i(t) \ge g_j(t)= \max_q g_q(t). 
  \]

The $C(i,j)$ are $18$ cones covering $\TM$, and, clearly, the intersection of any two of them is a face of each.
Therefore the cones $C(i,j)$, together with all their faces, are the cells of a complex of (convex rational polyhedral) cones subdividing $\TM$.
This complex is the chamber complex  of the $\SU{3}$ triple multiplicities  (see Section \ref{cones}).
  We denote this complex $\chK$.

\subsection{All 18 cones $C(i,j)$ are full-dimensional}\label{sec:fullD}
We prove this now.

For each of the cones $C(i,j)$, let $\upcone(i,j)$ be its inverse image under the projection $p: \coneBZ \rightarrow \RR^6$ that sends $BZ(t,x)$ (the BZ triangle defined in Figure \ref{BZtx}) to $t$.
Then $\upcone(i,j)$ is the set of points $BZ(t,x)$ fulfilling: ``$f_i(t) = \min_p f_p(t)$; $g_j(t) =\max_q g_q(t)$; and $f_i(t) \ge g_j(t)$''.
Rewrite the condition ``$f_i(t) = \min_p f_p(t)$'' as ``$f_i(t)-x = \min_p (f_p(t)-x)$''.
In terms of the coordinates $y_i$ and $z_j$ of the space~$\spanBZ$, this is ``$y_i = \min_p y_p$''.
Likewise, rewrite ``$g_j(t) =\max_q g_q(t)$'' as ``$x-g_j(t) =\min_q (x-g_q(t))$''; this is ``$z_j = \min_q z_q$''.
Finally, rewrite ``$f_i(t) \ge g_j(t)$'' as ``$f_i(t)-x \ge - (x-g_j(t))$'', which is ``$y_i \ge -z_j$'', a condition that always holds in $\coneBZ$.
As a conclusion, $\upcone(i,j)$ is the set of points of $\coneBZ$ that fulfill:'' $y_i = \min_p y_p$ and $z_j = \min_q z_q$''.
Now, the group $\upgroup$, that permutes the coordinates of the BZ triangles,  acts transitively on the pairs of coordinates $(y_i, z_j)$.
As a consequence, $\upgroup$ permutes transitively the cones $\upcone(i,j)$.
Accordingly, $\Gl$ permutes transitively the cones $C(i,j)$. The cones $C(i,j)$ are thus either all full-dimensional, or all degenerated.
But these cones cannot be all degenerated, because their union $\TM$ is full-dimensional.
We conclude that the cones $C(i,j)$ are all full-dimensional.

\subsection{All cells of the chamber complex are simplicial}\label{sec:simplicial}
This is what we prove now. Since any cell of the chamber complex is a face of some chamber, and any face of any simplicial cone is simplicial, it is enough to prove that all chambers are simplicial.

Consider a chamber $C(i,j)$. It is defined by the conditions: $f_i=\min_p f_p$; $g_j=\min_q g_q$; $f_i \ge g_j$. These conditions translate into the system of inequalities:
\[
\begin{cases}
  \forall p \neq i, \quad f_p \ge f_i  \quad \text{($2$ inequalities)},\\
  \forall q \neq j, \quad g_p \le g_j \quad \text{($5$ inequalities)},\\
  f_i \ge g_q \quad \text{($1$ inequality)}.
\end{cases}  
\]
Two of these inequalities are actually consequences of the other six. Indeed, from the equations $z_1 - z_4 = z_3-z_6 = z_5-z_2$ that hold on the linear span of the BZ triangles, and the equations $z_i=x-g_i(t)$,  follow the relations:
\begin{equation}\label{rel g}
g_1(t)-g_4(t) = g_5(t)-g_2(t) = g_3(t) - g_6(t).
\end{equation}
Using these relations, one gets that, for all $j$, 
\[
g_{j+1}-g_j = g_{j+1}-g_{j+4} + g_{j+4} - g_j = g_{j+3}-g_{j} + g_{j+4}-g_{j},
\]
where the indices are considered modulo $6$ (e.g. if $j=3$ then $j+4=1$).
Therefore $g_{j+1} \le g_j$ follows from $g_{j+3} \le g_j$ and $g_{j+4} \le g_j$.
Similarly, $g_{j-1} \le g_j$ follows from $g_{j+3} \le g_j$ and $g_{j+2} \le g_j$.
As a consequence, $C(i,j)$ is defined by the smaller system:
\begin{equation}\label{system Cij}
\begin{cases}
  f_{i+1} \ge f_i,\\
  f_{i+2} \ge f_i,\\
  g_{j+2} \le g_j, \\
  g_{j+3} \le g_j, \\
  g_{j+4} \le g_j,\\
  f_i \ge g_j.
  \end{cases}  
\end{equation}
Since $C(i,j)$ is a $6$-dimensional pointed cone, defined by a system of only $6$ linear inequalities, it is simplicial and all 6 inequalities are essential, i.e. the corresponding equations define the facets of $C(i,j)$.

\subsection{Combinatorial encoding of the cells}\label{sec:cells}

Let
 \[
 \Omega=\{C_1, C_2, C_3, D_1, D_3, D_5, \ltriangle, \rtriangle, \star\}
 \]
be the set of nine conditions defined in Table \ref{table:conditions}.

\begin{table}[ht]
  \centering
     \caption{The conditions whose conjunctions define the cells of the chamber complex $\chK$.}\label{table:conditions}

   \begin{tabular}{l@{\;:\;}l}
     \toprule
     $C_1$        & $\min_p f_p(t)=f_1(t)$ \\
     $C_2$        & $\min_p f_p(t)=f_2(t)$ \\
     $C_3$        & $\min_p f_p(t)=f_3(t)$ \\
     $D_1$        &  $\max_q g_q(t)=\max\left(g_1(t), g_{4}(t)\right)$\\
     $D_3$        & $\max_q g_q(t)=\max\left(g_3(t), g_{6}(t)\right)$\\
     $D_5$        & $\max_q g_q(t)=\max\left(g_5(t), g_{2}(t)\right)$\\
     $\ltriangle$ &  $\max_q g_q(t) = \max\left(g_1(t), g_3(t), g_5(t)\right)$\\
     $\rtriangle$ & $\max_q g_q(t) = \max \left(g_2(t), g_4(t), g_6(t)\right)$\\
     $\star$      & $\min_p f_p(t) = \max_q g_q(t)$\\
     \bottomrule
   \end{tabular}
   
 \end{table}

As a subset of $\TM$ (i.e. assuming that all inequalities $f_p(t) \geq g_q(t)$ are fulfilled), the chamber $C(i,j)$ is defined by
\begin{center}
  ``$\min_p f_p(t) = f_i(t)$ and $\max_q g_q(t) = g_j(t)$''.
 \end{center}
The first of the two conditions in this conjunction is exactly $C_i$. The second one is easily seen to be equivalent to
``$D_{j'}$ and $\omega_j$'', where $j'=j$ and $\omega_j=\rtriangle$ if $j$ is odd, and  $j'=j+3$ and $\omega_j=\ltriangle$ if $j$ is even.

Therefore, the chamber $C(i,j)$ is defined in $\TM$ by ``$C_i$ and $D_{j'}$ and $\omega_j$''.

Because $C(i, j)$  is simplicial, there is a one-to-one correspondence between the faces $\sigma$  of $C(i,j)$  and the sets $X$ of equations obtained from the defining inequalities \eqref{system Cij}, which are:
\begin{equation}\label{facet equations}
  f_{i+1} = f_i,\;
  f_{i+2} = f_i,\;
  g_{j+2} = g_j,\;
  g_{j+3} = g_j,\;
  g_{j+4} = g_j,\;
  f_i = g_j.
\end{equation}
In this correspondence, $\sigma$ corresponds to $X$ when $\sigma$ is defined, as a subset of $C(i,j)$, by the system of all equations in $X$;  and then $X$ is the set of all equations from the list that hold everywhere on $\sigma$.

The equations in \eqref{facet equations} can be restated in terms of conditions in $\Omega \setminus \{C_i, D_{j'}, \omega_j\}$. Indeed, the first two equations in \eqref{facet equations} are $C_{i+1}$ and $C_{i+2}$ (taking into account the assumption $f_i=\min_p f_p$). The last one is $\star$. One checks easily that the other three conditions ($g_{j+2} = g_j$, $g_{j+3} = g_j$,  $g_{j+4} = g_j$) are equivalent to $D_{j'+2}$, $\omega'_j$ and $D_{j'+4}$ respectively, where $\omega'_j$ is the element of $\{\ltriangle, \rtriangle\}$ that is not $\omega_j$. The above one-to-one correspondence becomes a one-to-one correspondence between the faces $\sigma$ of $C(i,j)$ and the subsets $X$ of  $\Omega \setminus \{C_i, D_{j'}, \omega_j\}$.

Note that, when $\sigma$ corresponds to $X$, then $\sigma$ is defined \emph{in $\TM$} by $X \cup \{C_i, D_{j'}, \omega_j\}$; and  $X \cup \{C_i, D_{j'}, \omega_j\}$ is the set of all conditions in $\Omega$ that hold on $\sigma$. We define two maps $\Psi$ and $\sol$: for any cell $\sigma$ of the chamber complex, $\Psi_{\sigma}$ is the set of all conditions in $\Omega$ that hold everywhere in $\sigma$;   for any subset $X\subset \Omega$, $\sol X$ is the set of points $t$ in $\TM$  that fulfill all conditions in $X$. From what precedes, we get that for any chamber $C(i,j)$, $\Psi$ and $\sol$ induce bijections, inverse of each other,  between the set of all faces of $C(i,j)$ and the subsets $X$ of $\Omega$ containing each of $C_i$, $D_{j'}$ and $\omega_j$.

  Together with the fact that each cell of $\chK$ is the face of some chamber, 
  this is enough to deduce that $\Psi$ and $\sol$ induce bijections, inverse of each other, between the set of all faces of the chamber complex, and the set $S$ of all subsets of $\Omega$ meeting each of $\{C_1, C_2, C_3\}$, $\{D_1, D_3, D_5\}$ and $\{\ltriangle, \rtriangle\}$.

  Rather than dealing with $\Psi$, which reverses inclusions, let us introduce, for any cell $\sigma$, the complement $\overline{\Psi}(\sigma)$ of $\Psi(\sigma)$ in $\Omega$ (the set of all conditions that fail to hold everywhere on $\sigma$).
  Then $\overline{\Psi}$ is an inclusion--preserving bijection from the set of all cells of $\chK$, to the set $\overline{S}$ of all parts of $\Omega$ that contain none of $\{C_1, C_2, C_3\}$, $\{D_1, D_3, D_5\}$ and $\{\ltriangle,\rtriangle\}$.
  The minimum element of $\overline{S}$ is the empty set, that corresponds to the $\{0\}$ cell.
  The minimal non-empty elements of $\overline{S}$ are the one-element subsets from $\Omega$, naturally in bijection with $\Omega$; they correspond to the rays of the chamber complex. Since each cell is simplicial, its dimension is the number of rays it contains.
  Therefore, the dimension of any cell $\sigma$ is the cardinality of $\overline{\Psi}(\sigma)$.

Given any set $E$, denote by $\PP{\, E}$ (resp. $\PPs{\, E}$) the set of all subsets $E$ (resp.  of all proper subsets of $E$, i.e. all subsets of $E$ distinct from $E$). The ranked poset of the cells of the chamber complex (ordered by inclusion) is isomorphic to the poset of the elements of $\overline{S}$, which is itself isomorphic to
\[
 \PPs{\{C_1, C_2, C_3\}} \times \PPs{\{D_1, D_3, D_5\}} \times \PPs{\{\ltriangle, \rtriangle\}} \times \PP{\{\star\}}.
 \]
 The generating series of a ranked poset is the polynomial $\sum_i a_i q^i$ where $q$ is a variable and $a_i$ is the number of elements of rank $i$.
 When $N$ is the cardinality of $E$, the generating series of $\PP{\, E}$ is $(1+q)^N$, and the generating series for $\PPs{\, E}$ is $(1+`q)^N-q^N$, which gives $1+3q+3q^2$ for $N=3$, and $(1+2q)$ when $N=2$.
 
Therefore,  $\overline{S}$ is a ranked poset with rank generating series
\[
(1+3q+3q^2)^2 (1+2 q) (1+q).
\]
This expands as 
\[
1 + 9 q + 35 q^2 + 75 q^3 + 93 q^4 + 63 q^5 + 18 q^6.
\]
The coefficients in this expansion are thus the numbers of faces of each dimension in the chamber complex of $SU(3)$ (``$f$-vector'').
In particular, we recover that the chamber complex has 18 chambers, and observe that it has nine rays. Eight of them have already been obtained, as the eight rays of the cone $\TM$. The group $G$ of linear symmetries of the chamber complex  permutes the nine rays of the chamber complex. On the other hand, we know that $G$ stabilizes the eight rays of the cone $\TM$. Therefore, $G$ fixes the ninth ray. In particular, this ray must be fixed by the six permutations of $(\ell;m;n)$ and by the duality involution. It follows that this ninth ray is generated by $(11|11|11)$. We denote the generator by $\ray{\star}$ and refer to its ray as the \emph{internal ray}, since it is the only one not on the border of $\TM$.

Each ray generator must fulfill all conditions from $\Omega$ but one.  One checks (using the expressions in coordinates in Table \ref{fundamental BZT}) that, for each $X \in \Omega$, the condition not fulfilled by the ray generator $\ray{X}$ is $X$. This justifies a posteriori the coincidence of notations.

Let us observe finally that the nine ray generators obtained here are related by
  \begin{equation}\label{rel between rays}
    \ray{C_1}+\ray{C_2}+\ray{C_3}
    =\ray{D_1}+\ray{D_3}+\ray{D_5}
    =\ray{\ltriangle}+\ray{\rtriangle}
    =\ray{\star}.
  \end{equation}
  Indeed, these relations are those in \eqref{rel between R}, except for
  $\ray{\ltriangle}+\ray{\rtriangle}=\ray{\star}$, which is, in coordinates,
  \[
  (01|01|01) + (10|10|10) = (11|11|11).
  \]

  Since these are $9$ vectors spanning a $6$--dimensional vector space, and since \eqref{rel between rays} already gives $3$ independent relations, there are no more independent relations.

\subsection{Cell diagrams and actions of the groups on the cells.}

There is a convenient pictorial way to describe the action of~$G$ and its subgroups on the cells. Draw a regular hexagon with vertices labelled with $\ray{C_1}$,  $\ray{D_1}$, $\ray{C_2}$, $\ray{D_3}$, $\ray{C_3}$ and $\ray{D_5}$, in this order.
Draw also the triangle whose vertices are the midpoints of the sides $\ray{C_2}\ray{D_1}$, $\ray{C_3}\ray{D_3}$ and $\ray{D_5} \ray{C_1}$ of the hexagon (``left--pointing triangle'') and the triangle whose vertices are the midpoints of the three remaining sides (``right--pointing triangle'').
Label accordingly these two triangles with $\ray{\ltriangle}$ and $\ray{\rtriangle}$.

\begin{figure}[ht]
  \centering
  \begin{tikzpicture}[scale=1.5]
  \node(0,0) (E)  [b, label={[label distance=-2mm] above right:$\ray{\star}$}] {};
  \path (E) -- ++(-30:1cm) node (C1) [b, label=below right:$\ray{C_1}$] {};
  \path (E) -- ++(+30:1cm) node (D1) [b,label=above right:$\ray{D_1}$] {};
  \path (E) -- ++(+90:1cm) node (C2) [b,label=above:$\ray{C_2}$] {};
  \path (E) -- ++(+150:1cm) node (D3) [b,label=above left:$\ray{D_3}$] {};
  \path (E) -- ++(+210:1cm) node (C3) [b,label=below left:$\ray{C_3}$] {};
  \path (E) -- ++(+270:1cm) node (D5) [b,label=below:$\ray{D_5}$] {};
  \draw (C1)--(D1)--(C2)--(D3)--(C3)--(D5)--(C1) -- cycle;
  \righttdiag
  \lefttdiag

  %\draw [->] (-.5,0) arc (180:0:0.5cm);
\end{tikzpicture}
\caption{The items associated to the rays of the chamber complex.}\label{ray hexagon}
\end{figure}
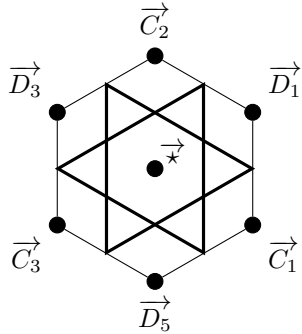

Label finally the center of this hexagon with  $\ray{\star}$,

We call the six vertices of the hexagon, its center and the two triangles the \emph{items of the cell diagrams}.
See Figure~\ref{ray hexagon}.

Given any cell of the chamber complex, its \emph{diagram} is obtained by superposing the items of the ray it contains. Table \ref{table rays} displays the diagrams of the nine rays, and Figure \ref{examples of diagrams} shows the diagrams of some other cells.

The missing items in the diagram of a cell $\sigma$ provide the conditions defining $\sigma$, i.e. the set $\Psi(\sigma)$. For instance, the diagram of the chamber $C(1,1)$ has all items except those of $C_1$, $D_1$ and $\ltriangle$ (see Figure \ref{examples of diagrams}). This corresponds to the fact that $C(1,1)$ is defined by the conditions $(C_1)$, $(D_1)$ and $(\ltriangle)$.

\begin{table}[ht]
  \centering
    \caption{The diagrams of the nine rays of the chamber complex, with the coordinates $(\ell_1,\ell_2|m_1,m_2|n_1, n_2)$ of their minimal generator.
    %  , and their name in Rassart's paper \cite{Rassart}
  }\label{table rays}
  \[
  \begin{array}{c@{\qquad}c@{\qquad}c}
    \toprule
    \begin{tikzpicture}[C1/.style={rb}, C2/.style={rb}, C3/.style={rb},D1/.style={rb}, D3/.style={b}, D5/.style={rb},E/.style={rb}] \fdiagram\end{tikzpicture}
      &
          \begin{tikzpicture}[C1/.style={rb}, C2/.style={rb}, C3/.style={rb},D1/.style={rb}, D3/.style={rb}, D5/.style={b},E/.style={rb}] \fdiagram\end{tikzpicture}
      &
            \begin{tikzpicture}[C1/.style={rb}, C2/.style={rb}, C3/.style={rb},D1/.style={b}, D3/.style={rb}, D5/.style={rb},E/.style={rb}] \fdiagram\end{tikzpicture}\\
                   \ray{D_3}=(00|10|01)         &  \ray{D_5}=(01|00|10) & \ray{D_1}=(10|01|00)\\[5mm]
%        (Rassart: $g_2$)      & (Rassart: $e_1$)        &  (Rassart: $d_2$) \\[5mm]
            \begin{tikzpicture}[C1/.style={b}, C2/.style={rb}, C3/.style={rb},D1/.style={rb}, D3/.style={rb}, D5/.style={rb},E/.style={rb}] \fdiagram\end{tikzpicture}
                &
            \begin{tikzpicture}[C1/.style={rb}, C2/.style={b}, C3/.style={rb},D1/.style={rb}, D3/.style={rb}, D5/.style={rb},E/.style={rb}] \fdiagram\end{tikzpicture}
            &
          \begin{tikzpicture}[C1/.style={rb}, C2/.style={rb}, C3/.style={b},D1/.style={rb}, D3/.style={rb}, D5/.style={rb},E/.style={rb}] \fdiagram\end{tikzpicture}
            \\
     \ray{C_1}=(00|01|10)   &   \ray{C_2}=(10|00|01) & \ray{C_3}=(01|10|00) \\[5mm]
%            (Rassart: $e_2$) &  (Rassart: $g_1$)  &  (Rassart: $d_1$)  \\[5mm]
                      \begin{tikzpicture}[C1/.style={rb}, C2/.style={rb}, C3/.style={rb},D1/.style={rb}, D3/.style={rb}, D5/.style={rb},E/.style={rb}] \fdiagram\righttdiag\end{tikzpicture}
                        &
                        \begin{tikzpicture}[C1/.style={rb}, C2/.style={rb}, C3/.style={rb},D1/.style={rb}, D3/.style={rb}, D5/.style={rb},E/.style={b}] \fdiagram\end{tikzpicture}
                        &
          \begin{tikzpicture}[C1/.style={rb}, C2/.style={rb}, C3/.style={rb},D1/.style={rb}, D3/.style={rb}, D5/.style={rb},E/.style={rb}] \fdiagram\lefttdiag\end{tikzpicture}
      \\
    \ray{\rtriangle}=(01|01|01) & \ray{\star}=(11|11|11) & \ray{\ltriangle}=(10|10|10) \\
    %   (Rassart: $c$)  &  (Rassart: $b$)  & (Rassart: $f$)
    \bottomrule
  \end{array}
  \]
  \end{table}

\begin{figure}[ht] \label{fig:chamber-complex-elements}
  \centering
  \[
  \begin{array}{c@{\qquad}c@{\qquad}c@{\qquad}c}
     \begin{tikzpicture}[C1/.style={rb}, C2/.style={b}, C3/.style={b},D1/.style={rb}, D3/.style={b}, D5/.style={b},E/.style={b}] \fdiagram\righttdiag\end{tikzpicture}
    &\begin{tikzpicture}[C1/.style={rb}, C2/.style={b}, C3/.style={b},D1/.style={rb}, D3/.style={b}, D5/.style={b},E/.style={rb}] \fdiagram\righttdiag\end{tikzpicture}
    &\begin{tikzpicture}[C1/.style={rb}, C2/.style={rb}, C3/.style={b},D1/.style={rb}, D3/.style={b}, D5/.style={rb},E/.style={rb}] \fdiagram\end{tikzpicture}
    &\begin{tikzpicture}[C1/.style={rb}, C2/.style={rb}, C3/.style={rb},D1/.style={rb}, D3/.style={rb}, D5/.style={rb},E/.style={rb}] \fdiagram\end{tikzpicture}
  \end{array}
  \]
      \caption{Diagrams of some cells of the chamber complex. From left to right: the diagram of the chamber $C(1,1)$, defined by the conditions $(C_1)$, $(D_1)$ and $\ltriangle$; the diagram of the external facet of $C(1,1)$, defined in $C(1,1)$ by $f_1=g_1$; the diagram of the $2$--dimensional face of $C(1,1)$ defined by the equations $f_1=f_2=g_1=g_4=g_5$; the diagram of the $\{0\}$ cell.}\label{examples of diagrams}
  \end{figure}

The action of the group $G$ and its subgroups on the cells is easily read from the diagrams.

Consider first the subgroup $\Glg$: it acts on the items as the group of symmetries of the triangle $\ray{C_1} \ray{C_2} \ray{C_3}$ (since $s_1$ acts as the reflection with respect to the axis $\ray{D_1} \ray{C_3}$, and $s_2$ as the reflection with respect to the axis $\ray{D_3} \ray{C_1}$.)

The duality involution acts on the items as the central symmetry (with respect to the center of the hexagon). Therefore, $\Gg$ (being generated by $\Glg$ and the duality involution) acts on the items as the group of the symmetries of the hexagon.

In what regards $\Gl$, recall that it permutes the rays as $$\S{\{\ray{C_1},\ray{C_2},\ray{C_3}\}} \times \S{\{\ray{D_1},\ray{D_3},\ray{D_5}\}} \times \S{\{\ray{\ltriangle}, \ray{\rtriangle}\}}.$$
Each factor of this decomposition acts as the group of symmetries of the triangle $C_1 C_2 C_3$ acting only on part of the diagram (leaving unaffected the other parts): 
 $\S{\{\ray{C_1},\ray{C_2},\ray{C_3}\}}$ on the triangle $\ray{C_1} \ray{C_2} \ray{C_3}$; $\S{\{\ray{D_1},\ray{D_3},\ray{D_5}\}}$ on the triangle $\ray{D_1} \ray{D_3} \ray{D_5}$; and $\S{\{\ray{\ltriangle}, \ray{\rtriangle}\}}$ on the left and right pointing triangles.

This describes the action of all elements of the group $G$, since $G$ is generated by $\Gl$ and the duality involution.

%%%%%%%%%%%%%%%%%%%%%%%%%%%%%%%%%%%%%%%%%%%%%%%%%%%%%%%%%%%%%%%%%%%%%%%%%%%%%%%
\subsubsection*{Example}

We consider the 3-dimensional cell represented by the following diagram and its orbits under the groups of symmetries:

\[
  \begin{tikzpicture}[scale=0.7,C1/.style={rb}, C2/.style={rb}, C3/.style={b},D1/.style={rb}, D3/.style={b}, D5/.style={rb},E/.style={rb}] \fdiagram\lefttdiag\end{tikzpicture}
\]

The elements of its orbit under $\Glg$, obtained by applying alternatively $s_1$ and $s_2$, are the cells with the following diagrams:
\[
\begin{array}{cccccc}
  \begin{tikzpicture}[scale=0.7,
      C1/.style={rb}, C2/.style={rb}, C3/.style={b},
      D1/.style={rb}, D3/.style={b}, D5/.style={rb},
      E/.style={rb}]
    \fdiagram\lefttdiag
  \end{tikzpicture}
  &
\begin{tikzpicture}[scale=0.7,
    C1/.style={rb}, C2/.style={b}, C3/.style={rb},
    D1/.style={rb}, D3/.style={b}, D5/.style={rb},
    E/.style={rb}]
  \fdiagram\righttdiag
\end{tikzpicture}
&
\begin{tikzpicture}[scale=0.7,
    C1/.style={b}, C2/.style={rb}, C3/.style={rb},
    D1/.style={rb}, D3/.style={rb}, D5/.style={b},
    E/.style={rb}]
  \fdiagram\lefttdiag
\end{tikzpicture}
&
\begin{tikzpicture}[scale=0.7,
    C1/.style={b}, C2/.style={rb}, C3/.style={rb},
    D1/.style={b}, D3/.style={rb}, D5/.style={rb},
    E/.style={rb}]
  \fdiagram\righttdiag
\end{tikzpicture}
&
\begin{tikzpicture}[scale=0.7,
    C1/.style={rb}, C2/.style={b}, C3/.style={rb},
    D1/.style={b}, D3/.style={rb}, D5/.style={rb},
    E/.style={rb}]
  \fdiagram\lefttdiag
\end{tikzpicture}
&
\begin{tikzpicture}[scale=0.7,
    C1/.style={rb}, C2/.style={rb}, C3/.style={b},
    D1/.style={rb}, D3/.style={rb}, D5/.style={b},
    E/.style={rb}]
  \fdiagram\righttdiag
\end{tikzpicture}
\end{array}
\]
This set of diagrams is stabilized by the central symmetry (duality involution for the cells). Therefore, the  orbit under $\Gg$ coincide with this orbit under $\Glg$.

The orbit under $\Gl$ and the orbit under $G$ also coincide, and correspond to the following set of diagrams, obtained by applying the hexagon symmetries only to the $C$-vertex, only to the $D$-vertex, or only to the triangles:
\[
\begin{array}{cccccc}
  \begin{tikzpicture}[scale=0.7,
      C1/.style={b}, C2/.style={rb}, C3/.style={rb},
      D1/.style={b}, D3/.style={rb}, D5/.style={rb},
      E/.style={rb}]
    \fdiagram\lefttdiag
  \end{tikzpicture}
  &
  \begin{tikzpicture}[scale=0.7,
      C1/.style={rb}, C2/.style={b}, C3/.style={rb},
      D1/.style={b}, D3/.style={rb}, D5/.style={rb},
      E/.style={rb}]
    \fdiagram\lefttdiag
  \end{tikzpicture}
  &
  \begin{tikzpicture}[scale=0.7,
      C1/.style={rb}, C2/.style={rb}, C3/.style={b},
      D1/.style={b}, D3/.style={rb}, D5/.style={rb},
      E/.style={rb}]
    \fdiagram\lefttdiag
  \end{tikzpicture}
  &
    \begin{tikzpicture}[scale=0.7,
      C1/.style={b}, C2/.style={rb}, C3/.style={rb},
      D1/.style={b}, D3/.style={rb}, D5/.style={rb},
      E/.style={rb}]
    \fdiagram\righttdiag
  \end{tikzpicture}
  &
  \begin{tikzpicture}[scale=0.7,
      C1/.style={rb}, C2/.style={b}, C3/.style={rb},
      D1/.style={b}, D3/.style={rb}, D5/.style={rb},
      E/.style={rb}]
    \fdiagram\righttdiag
  \end{tikzpicture}
  &
  \begin{tikzpicture}[scale=0.7,
      C1/.style={rb}, C2/.style={rb}, C3/.style={b},
      D1/.style={b}, D3/.style={rb}, D5/.style={rb},
      E/.style={rb}]
    \fdiagram\righttdiag
  \end{tikzpicture}
  \\
  \begin{tikzpicture}[scale=0.7,
      C1/.style={b}, C2/.style={rb}, C3/.style={rb},
      D1/.style={rb}, D3/.style={b}, D5/.style={rb},
      E/.style={rb}]
    \fdiagram\lefttdiag
  \end{tikzpicture}
  &
  \begin{tikzpicture}[scale=0.7,
      C1/.style={rb}, C2/.style={b}, C3/.style={rb},
      D1/.style={rb}, D3/.style={b}, D5/.style={rb},
      E/.style={rb}]
    \fdiagram\lefttdiag
  \end{tikzpicture}
  &
  \begin{tikzpicture}[scale=0.7,
      C1/.style={rb}, C2/.style={rb}, C3/.style={b},
      D1/.style={rb}, D3/.style={b}, D5/.style={rb},
      E/.style={rb}]
    \fdiagram\lefttdiag
  \end{tikzpicture}
  &
    \begin{tikzpicture}[scale=0.7,
      C1/.style={b}, C2/.style={rb}, C3/.style={rb},
      D1/.style={rb}, D3/.style={b}, D5/.style={rb},
      E/.style={rb}]
    \fdiagram\righttdiag
  \end{tikzpicture}
  &
  \begin{tikzpicture}[scale=0.7,
      C1/.style={rb}, C2/.style={b}, C3/.style={rb},
      D1/.style={rb}, D3/.style={b}, D5/.style={rb},
      E/.style={rb}]
    \fdiagram\righttdiag
  \end{tikzpicture}
  &
  \begin{tikzpicture}[scale=0.7,
      C1/.style={rb}, C2/.style={rb}, C3/.style={b},
      D1/.style={rb}, D3/.style={b}, D5/.style={rb},
      E/.style={rb}]
    \fdiagram\righttdiag
  \end{tikzpicture}
  \\
  \begin{tikzpicture}[scale=0.7,
      C1/.style={b}, C2/.style={rb}, C3/.style={rb},
      D1/.style={rb}, D3/.style={rb}, D5/.style={b},
      E/.style={rb}]
    \fdiagram\lefttdiag
  \end{tikzpicture}
  &
  \begin{tikzpicture}[scale=0.7,
      C1/.style={rb}, C2/.style={b}, C3/.style={rb},
      D1/.style={rb}, D3/.style={rb}, D5/.style={b},
      E/.style={rb}]
    \fdiagram\lefttdiag
  \end{tikzpicture}
  &
  \begin{tikzpicture}[scale=0.7,
      C1/.style={rb}, C2/.style={rb}, C3/.style={b},
      D1/.style={rb}, D3/.style={rb}, D5/.style={b},
      E/.style={rb}]
    \fdiagram\lefttdiag
  \end{tikzpicture}
  &
    \begin{tikzpicture}[scale=0.7,
      C1/.style={b}, C2/.style={rb}, C3/.style={rb},
      D1/.style={rb}, D3/.style={rb}, D5/.style={b},
      E/.style={rb}]
    \fdiagram\righttdiag
  \end{tikzpicture}
  &
  \begin{tikzpicture}[scale=0.7,
      C1/.style={rb}, C2/.style={b}, C3/.style={rb},
      D1/.style={rb}, D3/.style={rb}, D5/.style={b},
      E/.style={rb}]
    \fdiagram\righttdiag
  \end{tikzpicture}
  &
  \begin{tikzpicture}[scale=0.7,
      C1/.style={rb}, C2/.style={rb}, C3/.style={b},
      D1/.style={rb}, D3/.style={rb}, D5/.style={b},
      E/.style={rb}]
    \fdiagram\righttdiag
  \end{tikzpicture}
\end{array}
\]

%%%%%%%%%%%%%%%%%%%%%%%%%%%%%%%%%%%%%%%%%%%%%%%%%%%%%%%%%%%%%%%%%%%%%%%%%%%%%%%

\section{Linear symmetries and chamber complex for the Littlewood--Richardson coefficients of $\GL{3}$}\label{sec:LR}

In this section, we translate the results obtained for the $\SU{3}$ triple multiplicities into the setting of the Littlewood-Richardson coefficients (tensor multiplicities for the general linear groups $\GL{k}$, here for $\GL{3}$).

\subsection{From $\SU{3}$ to $\GL{3}$}

For any weakly decreasing sequence of integers $\lambda=(\lambda_1, \lambda_2, \lambda_3)$, let $W_{\lambda}$ be the irreducible representation of $\GL{3}$ whose highest weight sends the diagonal matrix with entries $x_1, x_2, x_3$ to $\lambda_1 x_2 +\lambda_2 x_2 + \lambda_3 x_3$.   All irreducible finite-dimensional representations of $\GL{3}$ are of this form.
The representation $W_{\lambda}$ is \emph{polynomial} when $\lambda_3\ge 0$, and then the label $\lambda$ is an integer partition.

The irreducible representations $V_{\ell}$ of $\SU{3}$ are exactly the restrictions of the irreducible representations of $W_{\lambda}$ of $\GL{3}$.
Precisely, for $\ell=(\ell_1, \ell_2)$ and $\lambda=(\lambda_1, \lambda_2, \lambda_3)$,
$V_{\ell}$ is the restriction of $W_{\lambda}$ when
\begin{equation}\label{Dynkin in partition}
\begin{cases}
  \ell_1=\lambda_1-\lambda_2,\\
  \ell_2=\lambda_2-\lambda_3.
\end{cases}
\end{equation}  
Note that all $\GL{3}$ irreducible representations $W_{\lambda_1+k, \lambda_2+k, \lambda_3+k}$, for $k$ integer (positive or negative), restrict to the same representation of $\SU{3}$.

A necessary condition for the $\GL{3}$ tensor  multiplicity $\mult{W_\nu}{W_\lambda\otimes W_\mu}$  (multiplicity of $W_\lambda$ in the tensor product $W_{\mu} \otimes W_{\nu}$)  to be non-zero is that $|\lambda|+|\mu|=|\nu|$ (where $|\lambda|$ stands for the sum of the coordinates of $\lambda$).
If $\lambda, \mu$ and $\nu$ are partitions (i.e. all three representations are polynomial), this tensor multiplicity is called a  {\em Littlewood--Richardson coefficient} and denoted $\lr{\lambda}{\mu}{\nu}$.

The relation 
\[
\mult{W_{\nu}}{W_{\lambda} \otimes W_{\nu}} = \mult{V_{n^*}}{V_\ell \otimes V_m} = \tm{\ell}{m}{n}
\]
holds when $V_{\ell}$, $V_m$ and $V_{n^*}$ are the restrictions of $W_{\lambda}$, $W_{\mu}$ and $W_{\nu}$ respectively, that is when:
\begin{equation}\label{lmn in partitions}
  \left\lbrace\begin{array}{lll}
  \ell_1 = \lambda_1-\lambda_2,  & m_1 = \mu_1 - \mu_2, & n_2 = \nu_1-\nu_2,\\
  \ell_2 = \lambda_2-\lambda_3,  & m_2 = \mu_2 - \mu_3, & n_1 = \nu_2 - \nu_3.
  \end{array}\right.
\end{equation}  
(We prefer to denote as $V_{n^*}$ the restriction of $W_{\nu}$, rather than as $V_n$, to get indices as simple as possible for the triple multiplicities).

\subsection{The chamber complex for the Littlewood--Richardson coefficients}

We now describe the Littlewood-Richardson coefficients and, more generally, the $\GL{3}$ tensor multiplicities, in the geometric language of cones.

The support of the $\GL{3}$ tensor multiplicities (the set of triples of labels $(\lambda;\mu;\nu)$ such that the corresponding tensor multiplicity is nonzero) is a set of integer points in $\RR^3 \times \RR^3 \times \RR^3$. This support is contained in the 8-dimensional subspace $\lspace{LRC}$ defined by the linear equation $|\lambda| + |\mu| = |\nu|$.
Let $\lattice{LRC}=\ZZ^{9} \cap \lspace{LRC}$ be the lattice of all integer points in $\lspace{LRC}$.

Formulas \eqref{lmn in partitions} define a projection $\pi:\lspace{LRC} \rightarrow \lspace{TM}$.
The space $\lspace{LRC}$ decomposes as
\begin{equation}\label{direct sum}
\lspace{LRC} = \lspace{LRC}^0 \oplus \ker \pi 
\end{equation}
where $\lspace{LRC}^0$ be the 6-dimensional subspace defined by the equations $\lambda_3 = \mu_3 = 0$, and $\ker \pi$ is the kernel of the projection $\pi$. It is a 2-dimensional subspace, with basis 
$(\ulrray, \vlrray)$ where
\[
\ulrray=(1,1,1; 0,0,0; 1,1,1), \quad \vlrray=(0,0,0;1,1,1;1,1,1).
\]
The projection $\pi$ maps isomorphically $\lspace{LRC}^0$ onto $\lspace{TM}$; its inverse defines an embedding $\iota: \lspace{TM} \hookrightarrow \lspace{LRC}$ with image $\lspace{LRC}^0$. 
The  decomposition \eqref{direct sum} together with the isomorphism $\lspace{LRC}^0 \cong \lspace{TM}$, define an isomorphism
\begin{equation}\label{iso}
 \lspace{\TM} \times \RR^2 \cong \lspace{LRC}.
\end{equation}

Under this isomorphism: $\lattice{TM} \times \ZZ^2$ corresponds to $\lattice{LRC}$;
the cone  $\TM \times \RR^2$ corresponds to a cone  $\Horn$ (known as the the \emph{Horn cone});
and $\latTM \times \ZZ^2$ corresponds to $\latof{\Horn}$ (the set of lattice points of the Horn cone), and this is the support of the $\GL{3}$ tensor multiplicities.
The Horn cone is better known as the set of spectra of triples of hermitian matrices $(A, B, C)$ that fulfill $A+B=C$, see \cite{Fulton:Horn, Bhatia}.

Under the isomorphism \eqref{iso} we have also that:
the cone  $\TM \times (\RR^+)^2$ corresponds to a cone  $\Horn^+$ (known as the \emph{positive Horn cone}), intersection of $\Horn$ with the subset defined by $\lambda_3 \ge 0$ and $\mu_3 \ge 0$; and the set $\latTM \times (\ZZ^+)^2$ corresponds to $\latof{\Horn^+}$ (lattice points of the positive Horn cone). This is the support of the Littlewood--Richardson coefficients.

The polynomial formula  $1+L_{i,j}$, that holds for the triple multiplicities in a chamber $C(i,j)$ of the chamber complex $\chc{TM}$, gives (by plugging the expressions \eqref{lmn in partitions} of $\ell$, $m$, $n$ in terms of $\lambda$, $\mu$, $\nu$) a formula for the $\GL{3}$ tensor multiplicities that holds in the cone that corresponds to $C(i,j) \times \RR^2$. These cones are the chambers of a chamber complex subdividing $\Horn$. Each of its chamber is isomorphic to some $C(i,j) \times \RR^2$.

The same formulas apply, of course, to the Littlewood--Richardson coefficients.
These are defined only for $\lambda_3,\mu_3, \nu_3$ nonnegative.
The domains for these formulas correspond thus to the products  $C(i,j) \times (\RR^+)^2$ and are the chambers of a chamber complex $\chKLR$ subdividing $\Horn^+$. We will denote with $C^{\bullet}(i,j)$ the chamber obtained from $C(i,j)$.
Precisely, it is the cone generated by $\iota(C(i,j))$ and  the vectors $\ulrray$, $\vlrray$.

The rays of $\chKLR$ are generated by the vectors $\ulrray$, $\vlrray$ and the images under the embedding $\iota$ of the ray generators for $\chc{TM}$. These are given in Table \ref{lrrays}. 

\begin{table}[ht]
  \centering
    \caption{The ray generators $u=(\ell_1\ell_2|m_1m_2|n_1n_2)$ for $\chK$, and the corresponding ray generators $v=(\lambda_1\lambda_2\lambda_3|\mu_1\mu_2\mu_3|\nu_1\nu_2\nu_3)$ for $\chKLR$. They are determined from each other by the conditions: $u = \pi(v)$ with $v\in \lspace{LRC}^0$; and  $v=\iota(u)$.}\label{lrrays}

  \[
  \begin{array}{ccc}
    \toprule
    \text{ray generator $u$ of $\chK$} & \text{ray generator $v$ of $\chKLR$} & \text{Name of $v$ in \cite{Rassart}.} \\[3mm]
    \midrule
          \ray{D_1}=(10|01|00) &
      \lrray{D_1}=(100|110|111) &
      d_2 \\[3mm]
      \ray{D_3}=(00|10|01) &
      \lrray{D_3}=(000|100|100) &
      g_2 \\[3mm]
      \ray{D_5}=(01|00|10) &
      \lrray{D_5}=(110|000|110)  &
      e_1 \\[3mm]
      \ray{C_1}=(00|01|10) &
      \lrray{C_1}=(000|110|110) &
      e_2 \\[3mm]
      \ray{C_2}=(10|00|01) &
      \lrray{C_2}=(100|000|100) &
      g_1 \\[3mm]
      \ray{C_3}=(01|10|00) &
      \lrray{C_3}=(110|100|111) &
       d_1 \\[3mm]
       \ray{\rtriangle}=(01|01|01)&
       \lrray{\rtriangle}=(110|110|211)&
      c \\[3mm]
      \ray{\star}=(11|11|11) &
      \lrray{\star}=(210|210|321) &
      b \\[3mm]
      \ray{\ltriangle}=(10|10|10) &
      \lrray{\ltriangle}=(100|100|110) &
      f \\[3mm]
      &
      \ulrray=(111|000|111)&
      a_1 \\[3mm]
      &
      \vlrray=(000|111|111) &
      a_2
      \\
        \bottomrule
  \end{array}
  \]
\end{table}  

The chambers of $\chKLR$ are $8$--dimensional and generated by $8$ rays (the embeddings of the 6 rays of the corresponding chamber from $\chc{TM}$, plus $\ulrray$ and $\vlrray$).
As a consequence, they are all simplicial.
The cells of $\chKLR$ are in order--preserving bijection with the cartesian product of the set of all chambers of $\chK$ with  $\PP{\{\ulrray, \vlrray\}}$.
Its rank generating function is thus the rank generating function of $\chK$, multiplied with $(1+q)^2$, which gives
\[
(1+3q+3q^2)^2 (1+2 q) (1+q)^3.
\]

\subsection{Linear symmetries of the Littlewood-Richardson coefficients.}

Recall the decomposition \eqref{direct sum}: $\lspace{LRC} = \lspace{LRC}^0 \oplus \ker \pi$.
Among the $11$ ray generators of the chamber complex $\chKLR$ (Table \ref{lrrays}), the vectors $\ulrray$ and $\vlrray$ form a basis of $\ker \pi$, while the other nine generators belong to $\lspace{LRC}^0$. Therefore,  $\ulrray$ and $\vlrray$ are not involved in any linear relation between these $11$ vectors. As what regards the other nine generators, they are related by
  \[
  \lrray{C_1}+\lrray{C_2}+\lrray{C_3}=\lrray{D_1}+\lrray{D_3}+\lrray{D_5}=\lrray{\ltriangle}+\lrray{\rtriangle}=\lrray{\star},
  \]
  that come from \eqref{rel between rays} by applying the embedding of $\lspace{\TM}$ into $\lspace{LRC}$.

 As a consequence, the group of linear automorphisms of the ``Littlewood--Ri\-chard\-son fun\-ction'' $(\lambda;\mu;\nu) \mapsto \lr{\lambda}{\mu}{\nu}$ decomposes as the product $G_1 \times G_2$, where $G_1$ is the subgroup of the linear automorphisms fixing $\lspace{LRC}^0$, and $G_2$ is the subgroup of the automorphisms fixing $\ulrray$ and $\vlrray$. The subgroup $G_2$ is isomorphic to the group $G$ of all symmetries of the triple multiplicities.
  The subgroup $G_1$ is not trivial: it contains the involution swapping  $\ulrray$ and $\vlrray$. It is thus isomorphic to $\S{2}$.

As a  conclusion, the group of all linear symmetries of the function $(\lambda;\mu;\nu) \mapsto \lr{\lambda}{\mu}{\nu}$ has order $288$, and is isomorphic to $\S{2} \times G$, and thus to $\S{2}\times \S{2} \times \left(\S{3} \wr \S{2}\right)$.

The generator of $G_1$, swapping  $\ulrray$ and $\vlrray$,
maps  $(\lambda_1,\lambda_2,\lambda_3| \mu_1, \mu_2,\mu_3 | \nu_1, \nu_2, \nu_3)$ to
\[
(\lambda_1-\delta,\lambda_2-\delta,\lambda_3-\delta\;|\;
\mu_1+\delta, \mu_2+\delta,\mu_3+\delta
\;|\; \nu_1, \nu_2, \nu_3)
\]
where $\delta=\lambda_3-\mu_3$.

The existence of this involution is easily interpreted: the representation $W_{(1^k)}$ of $\GL{k}$ is the one dimensional representation where $g\in \GL{k}$ acts by multiplication by $\det(g)^k$; for any $\lambda$,
\[
W_{\lambda+(1^k)} \cong W_{\lambda} \otimes W_{(1^k)},
\]
and therefore, for any $\lambda$ and $\mu$,
\[
W_{\lambda+(1^k)} \otimes W_{\mu} \cong W_{\lambda} \otimes W_{(1^k)} \otimes W_{\mu} \cong W_{\lambda} \otimes  W_{\mu + (1^k)}.
\]
This gives that for any $\lambda$, $\mu$ and $\nu$,
\[
\lr{\lambda+(1^k)}{\mu}{\nu} = \lr{\lambda}{\mu+(1^k)}{\nu}.
\]

%%%%%%%%%%%%%%%%%%%%%%%%%%%%%%%%%%%%%%%%%%%%%%%%%%%%%%%%%%%%%%%%%%%%%%%%%%%%%%%
%% Cone projections
%%----------------------------------------------------------------------

\section{Triple multiplicities and Littlewood--Richardson coefficients as volumes}\label{sec:volume}

In this section, we observe that the linear part of the $\SU{3}$ triple multiplicity $\tm{\ell}{m}{n}$ for can be interpreted as the volume of a parallelotope (higher--dimensional parallelepiped), or, equivalently, as a determinant.

Remember that any chamber $C(i, j)$ of $\chK$ has $6$ rays: one of them is interior (not on the border of $\TM$), with generator $\ray{\star}$; the other $5$ are exterior.

We contend that when $t=(\ell;m;n)$ lies in the chamber $C(i,j)$, we have
  \begin{equation}\label{volume tm}
  c(t) = 1 + \Vol{\lattice{TM}}\left(\Pi(t_1, \ldots, t_5, t)\right)
  \end{equation}
  where $t_1$, \ldots, $t_5$ are the minimal generators of the $5$ exterior rays of $C(i,j)$, 
  $\Pi(t_1, \ldots, t_5, t)$ is the parallelotope generated by them and $t$ (the set of all combinations $x_1 t_1 + \ldots + x_5 t_5 + x_6 t$ such that all $x_i$ fulfill $0 \le x_i \le 1$);
  and $\Vol{\lattice{TM}}$ is the volume with respect to the lattice $\lattice{TM}$, i.e. the volume normalized in such a way that the fundamental domains of the lattice have volume $1$.

  Let us check this. On $C(i,j)$, we have $c(t)=1+L_{i,j}(t)$, with $L_{i,j}$ linear. We should prove that $L_{i,j}$ coincides with the volume of the parallelotope.

  On $C(i,j)$, both $L_{i,j}(t)$ and $\Vol{\lattice{TM}}\left(\Pi(t_1, \ldots, t_5, t)\right)$ are linear in $t$; both vanish on the external facet of $C(i,j)$ (generated by $t_1$, \ldots, $t_5$). It
  is enough to check that both evaluate equally at $\ray{\star}$ to conclude that they are equal.
  We have  $c(\ray{\star})=2$. Therefore $L_{i,j}(\ray{\star})= 1$.
  
  We will show that  $\Vol{\lattice{TM}}\left(\Pi(t_1, \ldots, t_5, t)\right)=1$ as well, that is that $\Pi(t_1,\ldots, t_5, \star)$ is a fundamental domain for $\Lambda$. This amounts in showing that the set $B= \{t_1, \ldots, t_5, \ray{\star}\}$ is a basis for $\lattice{}$.

  The $9$ minimal generators for the rays of $\chK$ (table \ref{table rays}) span the lattice $\lattice{TM}$ (Section~\ref{sec:all symmetries}). Between them hold the relations \eqref{rel between rays}.

The set $B$ contains all minimal generators for the rays of $\chK$, with the exception of:
  one of the $\ray{C_i}$,
  one of the $\ray{D_j}$,
  and one of $\ray{\ltriangle}$ and $\ray{\rtriangle}$.
  This follows from the fact that the chambers are the maximal elements of the poset of cells of $\chK$, and from the description of this poset given in Section \ref{sec:cells}

But then, by \eqref{rel between rays}, the set $B$ also spans $\lattice{TM}$. Since $B$ has $6$ elements and $\lattice{TM}$ has rank $6$, we conclude that $B$ is a basis of $\lattice{TM}$, which was what was to be demonstrated.

Formula \eqref{volume tm} can be written using a determinant.
Let $M(t_1, t_2, \ldots)$ be the matrix whose columns give the coordinates $(\ell_1, \ell_2;m_1, m_2; n_1, n_2)$ of the vectors $t_1, t_2, \ldots$
Then 
\[
\left|\det M(t_1, \ldots, t_5, t)\right|=\Vol{\ZZ^6}\left(\Pi(t_1, \ldots, t_5, t)\right).
\]
But because $\lattice{TM}$ has index $3$ in $\ZZ^6$ (since $\lattice{TM}$ is defined by \eqref{lattice equation}), there is the relation $\Vol{\lattice{TM}} = 3\; \Vol{\ZZ^6}$.
As a consequence, for $t=(\ell,m,n) \in C(i, j)$,
\[
\tm{\ell}{m}{n} = 1 + \frac{1}{3}\left|\det M(t_1, \ldots, t_5, t)\right|.
\]

Analogous formulas, proved with similar arguments, hold for the Littlewood--Richardson coefficients of $\GL{3}$.
Any chamber $C^{\bullet}(i,j)$ of $\chc{LRC}$ has as minimal generators for its rays $\lrray{\star}$, $\ulrray$, $\vlrray$ and 5 other vectors $v_1$, \ldots, $v_5$ from Table \ref{lrrays}. 
For any vectors $w_1$, $w_2$, \ldots in $\lspace{LRC}$ let $M^{\bullet}(w_1, w_2, \ldots)$ be the matrix whose columns give the coordinates $(\lambda_1, \lambda_2, \lambda_3; \mu_1, \mu_2, \mu_3; \nu_1, \nu_2)$ (the component $\nu_3$ is dropped)  of $w_1$, $w_2$, \ldots
Then, for any lattice point $w=(\lambda; \mu; \nu) \in C^{\bullet}(i, j)$, there is
  \begin{align*}
    \lr{\lambda}{\mu}{\nu}
    &=  1 + \Vol{\lattice{LRC}}\left(\Pi(v_1, \ldots, v_5,  \ulrray, \vlrray, w)\right)\\
    &=  1+ \left|\det M^{\bullet}(v_1, \ldots, v_5,  \ulrray, \vlrray, w)\right|.
  \end{align*}

The interpretation of \eqref{volume tm} (except for the constant $1$)
as the continuous volume of a parallelotope provides an interesting visual explanation for the following fact: if any of $\ell_1,\ell_2,m_1,m_2,n_1,n_2$ are $0$, then the corresponding triple multiplicity is $1$. Namely, whenever one of the Dynkin labels vanishes, we see that the parallelotope has volume $0$ since the dimension drops. The case when $\ell_2 = 0$ is given by Pieri's rule.

%%%%%%%%%%%%%%%%%%%%%%%%%%%%%%%%%%%%%%%%%%%%%%%%%%%%%%%%%%%%%%%%%%%%%%%%%%%%%%%

\section{Stability}\label{sec:stability}

\emph{Stability} in representation theory refers to the following property exhibited by functions $F$ defining structural constants in terms of their labels: for any fixed label $u$ and any ``well-chosen'' label $v$, the sequence with general term $F(u + k v)$ (depending on the integer $k$) is eventually constant.
This property was observed first in some instances by Murnaghan \cite{Murnaghan:1938} for Kronecker coefficients (tensor multiplicities for the symmetric groups), then in greater generality in \cite{Stembridge, SamSnowden}, and widely generalized for other families of representation--theoretic settings \cite{SamSnowden, ChurchFarb:Stability}.
For triple multiplicities and Littlewood-Richardson coefficients, the stability property is easily explained thanks to the combinatorial models.

The case under consideration in this paper, of the triples multiplicities for $\SU{3}$, provides an even simpler  ``toy example''.

\textsc{Claim:} Let $u \in \latTM$ such that $c(u)=1$. Then, for any $t \in \latTM$,  there exists a positive integer $k_0$ such that
\[
\forall k \ge k_0, \quad c(t+ k u) = 1 + \min_{i,j} \{L_{i,j}(t) : u \in C(i,j)\}
\]
and in particular $c(t + ku)$ is independent of $k$.

\begin{proof}
Set $L(\ell, m, n)=c(\ell, m, n)-1$. 
  
  The key is to observe that there exists some chamber $C(i_0, j_0)$ and some index $k_0$ such that $t+k u \in C(i_0, j_0)$  for all $k \ge k_0$.
  Indeed, ``$t+ku \in C(i_0, j_0)$ for $k$ big enough'' is equivalent to ``$u + \varepsilon t \in C(i_0, j_0)$ for $\varepsilon > 0$ small enough'', since $C(i_0, j_0)$ is defined by linear inequalities (set $\varepsilon = 1/k$).   
  Then $u$ lies in the same chamber $C(i_0, j_0)$.
  Since $c(u)=1$, we have $L(u)=0$. Since $u \in C(i_0, j_0)$, this yields $L_{i_0, j_0}(u)=0$.
  On the other hand, for all $k\ge k_0$, 
\[
L(t+ k u)= L_{i_0, j_0}(t+ku)=L_{i_0, j_0}(t) +k L_{i_0, j_0}(u) = L_{i_0, j_0}(t).
\]

Further, we note that since $t+ku \in C(i_0, j_0)$, we have that $L_{i_0, j_0}(t+ku)=\min_{i,j} L_{i,j}(t+ku)$.

Let $C(i,j)$ be any chamber. We have thus
\[
L_{i_0, j_0}(t) = L_{i_0, j_0}(t+ku) \le L_{i_0, j_0}(t+ku) = L_{i_0, j_0}(t) + k L_{i_0, j_0}(u).
\]
Consider now any the chamber $C(i,j)$ containing $u$.
Then $L_{i,j}(u)=L(u)=0$.
Therefore $L_{i_0, j_0}(t) \le L_{i,j}(t)$.

This shows that $L_{i_0, j_0}(t)$ is the minimum of all $L_{i,j}(t)$ over all chambers $C(i,j)$ containing $u$.
\end{proof}

The result not only tells us that the sequence stabilizes, but also yields its stable value:  this is $1+L_{i,j}(t)$, where $C(i,j)$ is any of the chambers containing $u+ \varepsilon t$ for small positive $\varepsilon$.

%%%%%%%%%%%%%%%%%%%%%%%%%%%%%%%%%%%%%%%%%%%%%%%%%%%%%%%%%%%%%%%%%%%%%%

\section{Final remarks}\label{sec:remarks}

\subsection{Symmetries that cannot be lifted}\label{no lifting}

We have seen in Section \ref{sec:liftable} that, among the $12$ general symmetries for the $\SU{3}$ triple multiplicities, only $6$ of them  can be lifted to symmetries of the cone of the BZ triangles.
This followed from the calculation of the full group of symmetries of $\latBZ$.
Let us give here instead a short argument of this fact.
We will show that $\ell \leftrightarrow m$ does not lift to a symmetry of $\latBZ$.

Consider the the eight vectors $\ray{C_i}$, $\ray{D_j}$, $\ray{\ltriangle}$ and $\ray{\rtriangle}$, all defined in Table \ref{fundamental BZT}.  Each is the projection of a unique BZ triangle $\upray{C_i}$, $\upray{D_j}$, $\upray{\ltriangle}$ and $\upray{\rtriangle}$ (the ``fundamental" BZ triangles), also defined in Table \ref{fundamental BZT}.

Remember the following relation:
\begin{equation}
  \upray{D_1} + \upray{D_3} + \upray{D_5}
  =
  \upray{\ltriangle} + \upray{\rtriangle}.
  \tag{\ref{rel between fundamental}}
\end{equation}
On the other hand, the involution $\ell \leftrightarrow m$ fixes $\ray{\ltriangle}$ and $\ray{\rtriangle}$, and swaps $\ray{D_1}$, $\ray{D_3}$ and $\ray{D_5}$ with
$\ray{C_3}$, $\ray{C_2}$ and $\ray{C_1}$ respectively. But, last, as one can check:
\begin{equation}\label{nonrel C3C2C1}
  \upray{C_3} + \upray{C_2} + \upray{C_1}
  \neq 
  \upray{\ltriangle} + \upray{\rtriangle}.
\end{equation}
A would--be lifting of $\ell \leftrightarrow m$ would fix $\upray{\ltriangle}$ and $ \upray{\rtriangle}$, and map $\upray{D_1}$,  $\upray{D_3}$, $\upray{D_5}$ to  $\upray{C_3}$,  $\upray{C_2}$ and $\upray{C_1}$.
Then applying this lifting to \eqref{rel between fundamental} would give in \eqref{nonrel C3C2C1} an equality. Therefore \eqref{rel between fundamental} together with  \eqref{nonrel C3C2C1} provide a clear obstruction for the existence of a lifting of  $\ell \leftrightarrow m$.

This argument adapts straightforwardly for the Littlewood--Richardson coefficients, using for instance the hive model (see \cite{PakVallejo:cones}):  replace in the above formulas each ray generator  $\ray{X}$ of $\chK$ with the corresponding ray generator $\lrray{X}$ of $\chKLR$ (see Table \ref{lrrays}),  and the unique BZ triangle $\upray{X}$ above $\ray{X}$, with the unique hive above $\lrray{X}$.

The impossibility of lifting $\lambda \leftrightarrow \mu$ in the $\GL{4}$ case was already pointed out in~\cite{PakVallejo:cones}. The above calculation shows this already happens for $\GL{3}$.  

\subsection{The $\SU{2}$ case}

In the $\SU{2}$ case, each of $\ell$, $m$ and $n$ has only one coordinates. In this case, it follows from the Pieri rule that $c(\ell,m,n)$ counts the integer solutions $x$ of
\[
2 x = \ell + m + n
\text{ and }  x \le \ell
\text{ and } x \le m 
\text{ and }  x \ge 0.
\]
As a consequence, $c(\ell, m, n)=1$ when 
\[
\ell + m + n \equiv 0 \mod 2
\text{ and }  \ell \le m + n  
\text{ and }  m \le \ell + n
\text{ and }  n \le m + n,
\]
and otherwise $c(\ell,m,n)=0$.

\subsection{About the $\SU{4}$ case}

In contrast to the $\SU{3}$ case where the chamber complex for the triple multiplicities has $18$ chambers, in the $\SU{4}$  case, we have calculated that there are $67,769$ chambers.
In this case, the group of linear symmetries has order $12$, i.e. there are only the $12$ general linear symmetries.
Indeed, any symmetry of the triple multiplicities is also a symmetry of the corresponding cone $\TM$, but the cone $\TM$ for $\SU{4}$ affords only the $12$ general linear symmetries.

\subsection{Non linear symmetries}

Coquereaux and Zuber \cite{CoquereauxZuber} have studied in detail an additional symmetry for the $\SU{3}$ tensor multiplicities, that is not linear but piecewise linear (see also \cite{PelletierRessayre, Grinberg} for partial extensions of this result for $\SU{k}$, $k > 3$).
It would be interesting to study systematically such symmetries.

\subsection{Formulas with minima and maxima}

\begin{figure}[ht]
  \centering
  \includegraphics[scale=0.7]{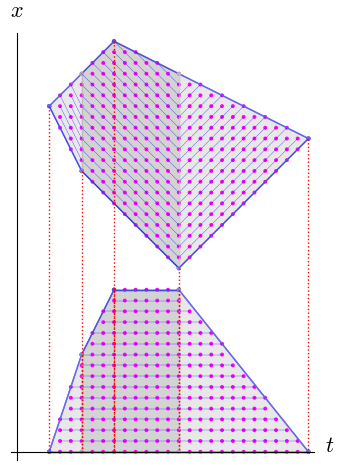}
  \caption{{\bf Top:} A transversal section of a cone defined by a system of inequalities $g_j(t) \le x \le f_i(t)$. {\bf Bottom:} The corresponding transversal version of its justification, defined by the inequalities $0 \le x \le f_i(t)-g_j(t)$.}\label{justification}
\end{figure}

The following two formulas for the triple multiplicities hold on $\TM$:
\begin{equation}\label{minmax formula}
c(t) = 1 + \min_i f_i(t) - \max_j g_j(t)
\end{equation}
and
\begin{equation}\label{min formula}
c(t) = 1 + \min_{i, j} \left(f_i(t)-g_j(t)\right).
\end{equation}
Both formulas appear in \cite[Formulas (16) and (17)]{CoquereauxZuber}.
Both interpret as counting lattice points in the fibers of a projection
$(x,t) \in \RR \times \RR^6 \mapsto t$ restricted to a cone.
The two formulas correspond to two different cones: for \eqref{minmax formula}, the cone is defined by the inequalities $f_i(t) \ge x$ and $x \ge g_j(t)$  (as in our presentation); for \eqref{min formula}, the cone is defined by the inequalities $0 \le x$ and $x \le f_i(t)-g_j(t)$, which is the region ``below the graph'' of the piecewise linear concave  function $\min_{i,j} (f_i(t)-g_j(t))$ . It is interesting to observe that the transformation mapping the first cone to the second, corresponds to a  simple  ``justification'' (to use the term from typography) of the fibers, i.e. pushing each fiber until its bottom is at level $0$, see Figure \ref{justification}. The parameters $\alpha$ used in \cite[section 3.7]{CoquereauxZuber} are precisely  the $x$--coordinates in the second description.

%% Important ingredients for the existence of the justification
%% 1) 1-dimensional
%% 2) the gi take integger values on intger points

%%%%%%%%%%%%%%%%%%%%%%%%%%%%%%%%%%%%%%%%%%%%%%%%%%%%%%%%%%%%%%%%%%%%%%
%% BIBLIOGRAPHY
%%----------------------------------------------------------------------

\bibliographystyle{hunsrt}
% or hunrst 
\bibliography{symmetries}

\begin{thebibliography}{10}

\bibitem{Elliott}
James~P. Elliott.
\newblock Collective motion in the nuclear shell model. {I.} {C}lassification
  schemes for states of mixed configurations.
\newblock {\em Proceedings of the Royal Society of London. Series A.
  Mathematical and Physical Sciences}, 245(1240):128--145, 1958.

\bibitem{Gellmann}
Murray {Gell-Mann}.
\newblock The eightfold way: a theory of strong interaction symmetry.
\newblock Technical Report TID-12608; CTSL-20, California Inst. of Tech.,
  Pasadena. Synchrotron Lab., 3 1961.

\bibitem{Neeman}
Yuval Ne'eman.
\newblock Derivation of strong interactions from a gauge invariance.
\newblock {\em Nuclear Physics}, 26(2):222--229, 1961.

\bibitem{Mandeltsveig}
Victor~B. Mandel{'}tsve\u{\i}g.
\newblock Irreducible representations of the {${\rm SU}_{3}$} group.
\newblock {\em Soviet Physics JETP}, 20:1237--1243, 1965.

\bibitem{ReckZeilingerBernstein}
Michael Reck, Anton Zeilinger, Herbert~J. Bernstein, and Philip Bertani.
\newblock Experimental realization of any discrete unitary operator.
\newblock {\em Phys. Rev. Lett.}, 73:58--61, Jul 1994.

\bibitem{Wesslen}
Maria S.~M. Wessl\'{e}n.
\newblock A geometric description of tensor product decompositions in
  {$\mathfrak{su}(3)$}.
\newblock {\em J. Math. Phys.}, 49(7):073506, 16, 2008.

\bibitem{Sanders_et_al}
Barry~C Sanders, Hubert de~Guise, D~J Rowe, and A~Mann.
\newblock Vector phase measurement in multipath quantum interferometry.
\newblock {\em Journal of Physics A: Mathematical and General},
  32(44):7791--7801, oct 1999.

\bibitem{Fulton}
William Fulton.
\newblock {\em {Y}oung tableaux}, volume~35 of {\em London Mathematical Society
  Student Texts}.
\newblock Cambridge University Press, Cambridge, 1997.

\bibitem{BerensteinZelevinsky}
Arkady~D. Berenstein and Andrei~V. Zelevinsky.
\newblock Triple multiplicities for {${\it sl}(r+1)$} and the spectrum of the
  exterior algebra of the adjoint representation.
\newblock {\em J. Algebraic Combin.}, 1(1):7--22, 1992.

\bibitem{BriandRosas}
Emmanuel Briand and Mercedes Rosas.
\newblock The 144 symmetries of the {L}ittlewood-{R}ichardson coefficients of
  ${SL_3}$, 2020, arxiv:2004.04995 [math.CO].

\bibitem{CrampeEtAl}
Nicolas Crampe, Lo\"{\i}c~Poulain d'Andecy, and Luc Vinet.
\newblock The missing label of $\mathfrak{su}_3$ and its symmetry.
\newblock {\em Communications in Mathematical Physics}, 2023.
\newblock Preliminary version (2021) at arXiv:2110.03521 [math.RT].

\bibitem{Rassart}
\'Etienne Rassart.
\newblock A polynomiality property for {L}ittlewood-{R}ichardson coefficients.
\newblock {\em J. Combin. Theory Ser. A}, 107(2):161--179, 2004.

\bibitem{PakVallejo:cones}
Igor Pak and Ernesto Vallejo.
\newblock Combinatorics and geometry of {L}ittlewood-{R}ichardson cones.
\newblock {\em European J. Combin.}, 26(6):995--1008, 2005.

\bibitem{CoquereauxZuber}
Robert Coquereaux and Jean-Bernard Zuber.
\newblock Conjugation properties of tensor product multiplicities.
\newblock {\em Journal of Physics A: Mathematical and Theoretical},
  47(45):455202, oct 2014.

\bibitem{DummitFoote}
David~S. Dummit and Richard~M. Foote.
\newblock {\em Abstract algebra}.
\newblock John Wiley \& Sons, Inc., Hoboken, NJ, third edition, 2004.

\bibitem{Fulton:Horn}
William Fulton.
\newblock Eigenvalues, invariant factors, highest weights, and {S}chubert
  calculus.
\newblock {\em Bull. Amer. Math. Soc. (N.S.)}, 37(3):209--249, 2000.

\bibitem{Bhatia}
Rajendra Bhatia.
\newblock Algebraic geometry solves an old matrix problem.
\newblock {\em Resonance}, 4:101--105, 1999.

\bibitem{Murnaghan:1938}
Francis~D. Murnaghan.
\newblock The {A}nalysis of the {K}ronecker {P}roduct of {I}rreducible
  {R}epresentations of the {S}ymmetric {G}roup.
\newblock {\em Amer. J. Math.}, 60(3):761--784, 1938.

\bibitem{Stembridge}
John Stembridge.
\newblock Generalized stability of {K}ronecker coefficients, August 2014.
\newblock Preprint num. 64 at
  \url{http://www.math.lsa.umich.edu/~jrs/papers.html}.

\bibitem{SamSnowden}
Steven~V. Sam and Andrew Snowden.
\newblock Stability patterns in representation theory.
\newblock {\em Forum Math. Sigma}, 3:Paper No. e11, 108, 2015.

\bibitem{ChurchFarb:Stability}
Thomas Church and Benson Farb.
\newblock Representation theory and homological stability.
\newblock {\em Adv. Math.}, 245:250--314, 2013.

\bibitem{PelletierRessayre}
Maxime Pelletier and Nicolas Ressayre.
\newblock Some unexpected properties of {L}ittlewood-{R}ichardson coefficients.
\newblock {\em Electronic Journal of Combinatorics}, 29(4):11, 2022.
\newblock Preliminary version (2020) at arxiv:2005.09877 [math.CO].

\bibitem{Grinberg}
Darij Grinberg.
\newblock The {P}elletier-{R}essayre hidden symmetry for
  {L}ittlewood-{R}ichardson coefficients.
\newblock {\em Combinatorial Theory}, 1:16, 2021.
\newblock Extended version at arxiv:2008.06128 [math:CO].

\end{thebibliography}

\end{document}